\input amssym.def
\input psfig
\input epsf

\let \blskip = \baselineskip
\parskip=1.2ex plus .2ex minus .1ex

\tabskip 20pt
\tolerance = 1000
\pretolerance = 50
\newcount\itemnum
\itemnum = 0
\overfullrule = 0pt

\def\title#1{\bigskip\centerline{\bigbigbf#1}}
\def\author#1{\bigskip\centerline{\bf #1}\smallskip}
\def\address#1{\centerline{\it#1}}
\def\abstract#1{\vskip1truecm{\narrower\noindent{\bf Abstract.} #1\bigskip}}

\def\sp{\bigskip}
\def\nosp{\vskip -\the\blskip plus 1pt minus 1pt}

\def\br{\hfil\break} 
\def\ti{\br \hglue \the \parindent}

\def\ce#1{\LP\centerline{#1}}

\def\skipit#1{}
\def\mdag{\raise 3pt\hbox{\dag}}

\def\XP{\par\noindent\hang}
\def\LP{\par\noindent}
\def\BP[#1]{\par\item{[#1]}}
\def\SH#1{\sp\vskip\parskip\leftline{\bigbf #1}\nobreak}

\def\TH#1{\sp\XP{\bf THEOREM\ \shead#1}}
\def\LM#1{\sp\XP{\bf LEMMA\ \shead#1}}
\def\DF#1{\sp\XP{\bf DEFINITION\ \shead#1}}

\def\PF{\LP{\bf Proof:\ }}
\def\NX{\advance\itemnum by 1 \sp\LP {\bf \shead \the\itemnum.\ }}
\def\qed{\null\nobreak\hfill\hbox{${\vrule width 5pt height 6pt}$}\par\sp}

\def\cart{\>\hbox{${\vcenter{\vbox{
    \hrule height 0.4pt\hbox{\vrule width 0.4pt height 4.5pt
    \kern4pt\vrule width 0.4pt}\hrule height 0.4pt}}}$}\>}
\def\bxmu{\>\hbox{${\vcenter{\vbox {
    \hrule height 0.4pt\hbox{\vrule width 0.4pt height 4pt
    \hskip -1.3pt\lower 1.8pt\hbox{$\times$}\negthinspace\vrule width 0.4pt}
    \hrule height 0.4pt}}}$}\>}

\def\lin#1{\hbox to #1true in{\hrulefill}}


\def\foot#1{\raise 6pt \hbox{#1} \kern -3pt}

\def\fig #1 #2 #3 #4 #5 {\sp \ce{ {\epsfbox[#1 #2 #3 #4]{figs/#5.ps}} }}

\def\gpic#1{#1 \sp\ce{\box\graph} \medskip} 


\def\JGT{{\it J.\ Graph Theory}}

\def\DM{{\it Discrete Math.{}}}
\def\DAM{{\it Discrete Appl.\ Math.{}}}

\def\ADM{{\it Annals Discr.\ Math.{}}}

\def\SIAD{{\it SIAM J.\ Algeb.\ Disc.\ Meth.{}}}
\def\SIAP{{\it SIAM J.\ Appl.\ Math.{}}}
\def\SIDM{{\it SIAM J.\ Discr.\ Math.{}}}

\def\CNum{{\it Congr.\ Numer.{}}}

\def\al{\alpha}	\def\be{\beta}		
    \def\lmb{\lambda}

\def\ZZ{{\Bbb Z}}  

\def\nul{\hbox{\O}}    
		
\def\esub{\subseteq}



\def\({\left(}	\def\){\right)}


\def\FR#1#2{{#1 \over #2}}

\def\CL#1{\left\lceil{#1}\right\rceil}

\def\SE#1#2#3{\sum_{#1 = #2} ^ {#3}}

\def\UE#1#2#3{\bigcup_{#1 = #2} ^ {#3}}

\def\SM#1#2{\sum_{#1\in#2}}

\def\UM#1#2{\bigcup_{#1\in#2}}

\def\VEC#1#2#3{#1_{#2},\ldots,#1_{#3}}

\def\st{\colon\;} 
\def\MAP#1#2#3{#1\colon\;#2\to#3}
\def\SET#1:#2{\{#1\colon\;#2\}}


		
\def\C#1{\left | #1 \right |}    

        
    \def\diam{ {\rm diam\;}}





\magnification=\magstep1
\vsize=9.0 true in
\hsize=6.5 true in
\headline={\hfil\ifnum\pageno=1\else\folio\fi\hfil}
\footline={\hfil\ifnum\pageno=1\folio\else\fi\hfil}

\parindent=20pt
\baselineskip=12pt
\parskip=.5ex  

\def\shead{ }

\font\bigbf = cmb10 scaled \magstep1

\font\bigbigbf = cmb10 scaled \magstep2


\title{BANDWIDTH AND DENSITY FOR BLOCK GRAPHS}
\author{Le Tu Quoc Hung and Maciej M. Sys\l o}
\address{Institute of Computer Science, University of Wroc\l aw}
\address{Przesmyckiego 20, 51-151 Wroc\l aw, Poland}
\author{Margaret L. Weaver}
\address{Mathematics Department, Eastern Illinois University}
\address{Charleston, IL 61920-3099}
\author{Douglas B. West\foot{\dag}}
\address{Mathematics Department, University of Illinois}
\address{Urbana, IL 61801-2975}
\vfootnote{}{\br
   \foot{\dag}Research supported in part by NSA/MSP Grant MDA904-93-H-3040.\br
   Running head: {BANDWIDTH AND BLOCK GRAPHS}\br
   AMS codes: 05C78\br
   Keywords: bandwidth, labeling, block graph, NP-complete\br
   Written November 1992; revised October 1996.
}
\abstract{The {\it bandwidth} of a graph $G$ is the minimum of the maximum
difference between adjacent labels when the vertices have distinct integer
labels.  We provide a polynomial algorithm to produce an optimal bandwidth
labeling for graphs in a special class of block graphs (graphs in which
every block is a clique), namely those where
deleting the vertices of degree one produces a path of cliques.  The result
is best possible in various ways.  Furthermore, for two classes of graphs
that are ``almost'' caterpillars, the bandwidth problem is NP-complete.}

\def\rep{representation}

\SH
{1. INTRODUCTION}
The {\it bandwidth problem} for a graph asks for a linear layout to minimize
stretching of edges (see [10] for a VLSI circuit layout application).
The {\it bandwidth} of an injection $\MAP f{V(G)}{\ZZ}$ is
$B(f) = \max_{uv\in E(G)} \C{f(u)-f(v)}$.  The {\it bandwidth} $B(G)$ of a
graph $G$ is $\min B(f)$ over all such injections; a numbering achieving the
minimum is {\it optimal}.  Surveys on bandwidth include [2] and [3].

Let $n(G) = \C{V(G)}$.  Every numbering of $G$ uses two labels differing by at
least $n(G)-1$, and the two corresponding vertices are connected by a path of
length at most $\diam G$; thus $B(G) \ge (n(G)-1)/\diam G$.  Considering all
subgraphs, the {\it local density} is 
$\be(G)=\max_{H\esub G} \CL{(n(H)-1)/\diam H}$.  Since every numbering of $G$
includes a numbering of each subgraph, $B(G)\ge\be(G)$ (see Chung [3]).  The
local density bound is optimal for cliques, stars, and trees of diameter 3
(``double stars'').  Sys\l o and Zak [17] and Miller [11] extended this to {\it
caterpillars}, the trees in which deleting the vertices of degree one produces a
path.  Their proofs construct optimal numberings in polynomial time.  This was
further extended by Assmann, Peck, Sys\l o, and Zak [1] to {\it 2-caterpillars}.
(A $k$-caterpillar is a tree formed from a path by growing edge-disjoint paths
of lengths at most $k$ from its vertices.  Commonly called ``caterpillars with
hairs of length at most $k$'', these are not caterpillars when $k\ge2$).
Among trees, bandwidth has also been computed for complete $k$-ary trees [15].

We extend the caterpillar result.  A graph is a {\it block graph} if every
block is a clique.  This name arises because a graph $G$ is the intersection
graph of the blocks of some graph if and only if every block of $G$ is a
clique [8].  A {\it block path} is a block graph with $k$ cutvertices
and $k+1$ blocks in which the cutvertices induce a path.  A {\it block
caterpillar} is a block graph in which deleting the leaves (1-valent vertices)
produces a block path.  Fig.~1 illustrates a block caterpillar; the ellipses
represent blocks that are cliques.  We provide an algorithm to construct
optimal numberings (with bandwidth $\be(G)$) for block caterpillars.  Note that
2-caterpillars are not generally block caterpillars.  We also demonstrate that
the local density bound need not be optimal for block graphs of diameter 3
or for trees of diameter 4.

\sp
\gpic{
\expandafter\ifx\csname graph\endcsname\relax \csname newbox\endcsname\graph\fi
\expandafter\ifx\csname graphtemp\endcsname\relax \csname newdimen\endcsname\graphtemp\fi
\setbox\graph=\vtop{\vskip 0pt\hbox{%
    \graphtemp=.5ex\advance\graphtemp by 1.069in
    \rlap{\kern 0.411in\lower\graphtemp\hbox to 0pt{\hss $\bullet$\hss}}%
    \graphtemp=.5ex\advance\graphtemp by 1.069in
    \rlap{\kern 0.740in\lower\graphtemp\hbox to 0pt{\hss $\bullet$\hss}}%
    \graphtemp=.5ex\advance\graphtemp by 1.069in
    \rlap{\kern 1.069in\lower\graphtemp\hbox to 0pt{\hss $\bullet$\hss}}%
    \graphtemp=.5ex\advance\graphtemp by 1.069in
    \rlap{\kern 1.398in\lower\graphtemp\hbox to 0pt{\hss $\bullet$\hss}}%
    \graphtemp=.5ex\advance\graphtemp by 1.069in
    \rlap{\kern 1.727in\lower\graphtemp\hbox to 0pt{\hss $\bullet$\hss}}%
    \graphtemp=.5ex\advance\graphtemp by 1.069in
    \rlap{\kern 2.056in\lower\graphtemp\hbox to 0pt{\hss $\bullet$\hss}}%
    \graphtemp=.5ex\advance\graphtemp by 1.069in
    \rlap{\kern 2.713in\lower\graphtemp\hbox to 0pt{\hss $\bullet$\hss}}%
    \graphtemp=.5ex\advance\graphtemp by 1.069in
    \rlap{\kern 3.371in\lower\graphtemp\hbox to 0pt{\hss $\bullet$\hss}}%
    \graphtemp=.5ex\advance\graphtemp by 0.411in
    \rlap{\kern 3.371in\lower\graphtemp\hbox to 0pt{\hss $\bullet$\hss}}%
    \graphtemp=.5ex\advance\graphtemp by 0.740in
    \rlap{\kern 0.082in\lower\graphtemp\hbox to 0pt{\hss $\bullet$\hss}}%
    \graphtemp=.5ex\advance\graphtemp by 0.082in
    \rlap{\kern 0.082in\lower\graphtemp\hbox to 0pt{\hss $\bullet$\hss}}%
    \graphtemp=.5ex\advance\graphtemp by 0.411in
    \rlap{\kern 0.411in\lower\graphtemp\hbox to 0pt{\hss $\bullet$\hss}}%
    \graphtemp=.5ex\advance\graphtemp by 0.411in
    \rlap{\kern 0.576in\lower\graphtemp\hbox to 0pt{\hss $\bullet$\hss}}%
    \graphtemp=.5ex\advance\graphtemp by 0.411in
    \rlap{\kern 1.069in\lower\graphtemp\hbox to 0pt{\hss $\bullet$\hss}}%
    \graphtemp=.5ex\advance\graphtemp by 0.411in
    \rlap{\kern 1.398in\lower\graphtemp\hbox to 0pt{\hss $\bullet$\hss}}%
    \graphtemp=.5ex\advance\graphtemp by 0.411in
    \rlap{\kern 1.891in\lower\graphtemp\hbox to 0pt{\hss $\bullet$\hss}}%
    \graphtemp=.5ex\advance\graphtemp by 0.411in
    \rlap{\kern 2.384in\lower\graphtemp\hbox to 0pt{\hss $\bullet$\hss}}%
    \graphtemp=.5ex\advance\graphtemp by 0.411in
    \rlap{\kern 2.713in\lower\graphtemp\hbox to 0pt{\hss $\bullet$\hss}}%
    \special{pn 8}%
    \special{pa 82 740}%
    \special{pa 411 411}%
    \special{pa 82 82}%
    \special{fp}%
    \special{pa 411 1069}%
    \special{pa 576 411}%
    \special{pa 740 1069}%
    \special{fp}%
    \special{pa 1727 1069}%
    \special{pa 1891 411}%
    \special{pa 2056 1069}%
    \special{fp}%
    \special{pa 2384 411}%
    \special{pa 2713 411}%
    \special{fp}%
    \special{pa 1069 1069}%
    \special{pa 1069 411}%
    \special{fp}%
    \special{pa 1398 1069}%
    \special{pa 1398 411}%
    \special{fp}%
    \special{pa 2713 1069}%
    \special{pa 2713 411}%
    \special{fp}%
    \special{pa 3371 1069}%
    \special{pa 3371 411}%
    \special{fp}%
    \special{ar 904 411 493 329 0 6.28319}%
    \special{ar 1891 411 493 329 0 6.28319}%
    \special{ar 3207 411 493 329 0 6.28319}%
    \hbox{\vrule depth1.151in width0pt height 0pt}%
    \kern 3.700in
  }%
}%
}
\vskip .5pc
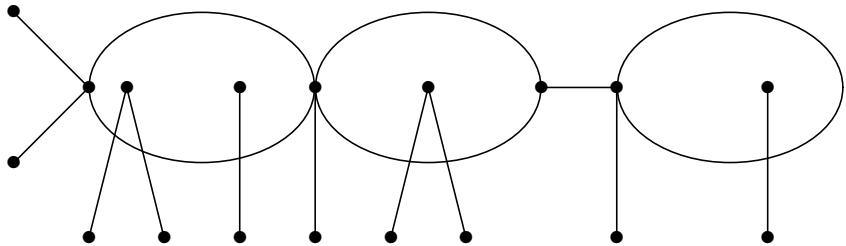
\ce{Fig. 1.  Sketch of a block caterpillar}

Computing bandwidth is NP-complete [14], even for trees with maximum degree 3
[4]; hence the interest in special classes.  Slightly enlarging the classes of
2-caterpillars or block caterpillars yields classes on which bandwidth is
NP-complete.  Monien [12] proved that bandwidth is NP-complete for
3-caterpillars, although he needs paths of length 3 only at one vertex of the
central path.  We prove NP-completeness for two additional classes.
One class consists of block graphs obtained from special block caterpillars
by adding paths of length two from one vertex of the central path.  The other
class consists of trees that are almost caterpillars; they have a path
containing all non-leaf vertices except one.

Because these trees are tolerance graphs, we conclude that bandwidth is
NP-complete for tolerance graphs, answering a question posed by Kleitman.
A graph is a {\it tolerance graph} if it is possible to assign each vertex $v$
an interval $I_v=[a_v,b_v]$ and a tolerance $t_v$ such that vertices $x,y$ are
adjacent if and only if $I_x\cap I_y$ has length at least $\min\{t_x,t_y\}$.
The class of tolerance graphs (introduced in [6] and [7]) contains the class
of interval graphs, on which there are polynomial time algorithms for bandwidth
[9,13,16].  (Interval graphs are the tolerance graphs representable using the same
tolerance for all vertices; more simply, they are the intersection graphs of
families of real intervals.)
\skipit{
Of particular interest among the trees with maximum degree 3 are complete binary
trees.  Numberings achieving $\be(G)$ were proposed for complete binary trees
and for complete $t$-ary trees ([2,4]), flaws were found [5], and correct
numberings appear in [15].
}

\SH
{2. EXAMPLES}
Before proving the main result, we exhibit examples where bandwidth does not 
equal local density.  Equality holds for all trees of diameter 3, which are
caterpillars, but this does not extend to block graphs with diameter 3 or to
trees with diameter 4.

\TH 1.
There are block graphs of diameter 3 for which the bandwidth exceeds the local
density bound.
\PF
Consider the block graph $H_k$ with four blocks illustrated in Fig. 2.  Three
of the blocks are disjoint cliques $X,Y,Z$ of order $k$.  The fourth consists
of $x\in X$, $y\in Y$, $z\in Z$, and one additional vertex $w$ not in the other
cliques.  The largest subgraphs of $H_k$ having diameter $d$ have $k,k+3,3k+1$
vertices for $d = 1,2,3$, respectively (if $k\ge 3$), so $\be(H_k)=k$
when $k\ge 3$.  (For $k=2$ the graph is a block caterpillar, and $\be(H_2)=3$).

Suppose that $B(H_k)=k$.  We may assume that the optimal labeling $f$ uses
labels $\{0,\ldots,3k\}$.  The distance between $f^{-1}(0)$ and $f^{-1}(3k)$
must be 3, so we may assume that $f^{-1}(0)\in X$ and $f^{-1}(3k)\in Z$.  Hence
$f(x)=k$ and $f(z)=2k$.  Since every vertex is within distance 2 of $w$, we have
$k\le f(w)\le2k$.  Since $\C{X-x}=k-1$ and $\C{Z-z}=k-1$, we must now have
distinct vertices in $Y$ with labels less than $k$ and greater than $2k$.  This
yields adjacent vertices whose labels differ by more than $k$.  (Note: $w$ is
needed in this construction, since $B(H_k-w)=k$.)
\qed

\sp
\gpic{
\expandafter\ifx\csname graph\endcsname\relax \csname newbox\endcsname\graph\fi
\expandafter\ifx\csname graphtemp\endcsname\relax \csname newdimen\endcsname\graphtemp\fi
\setbox\graph=\vtop{\vskip 0pt\hbox{%
    \graphtemp=.5ex\advance\graphtemp by 0.732in
    \rlap{\kern 1.000in\lower\graphtemp\hbox to 0pt{\hss $\bullet$\hss}}%
    \graphtemp=.5ex\advance\graphtemp by 1.281in
    \rlap{\kern 1.317in\lower\graphtemp\hbox to 0pt{\hss $\bullet$\hss}}%
    \graphtemp=.5ex\advance\graphtemp by 1.281in
    \rlap{\kern 0.683in\lower\graphtemp\hbox to 0pt{\hss $\bullet$\hss}}%
    \graphtemp=.5ex\advance\graphtemp by 0.732in
    \rlap{\kern 1.000in\lower\graphtemp\hbox to 0pt{\hss $\bullet$\hss}}%
    \graphtemp=.5ex\advance\graphtemp by 1.281in
    \rlap{\kern 1.317in\lower\graphtemp\hbox to 0pt{\hss $\bullet$\hss}}%
    \graphtemp=.5ex\advance\graphtemp by 1.281in
    \rlap{\kern 0.683in\lower\graphtemp\hbox to 0pt{\hss $\bullet$\hss}}%
    \graphtemp=.5ex\advance\graphtemp by 1.098in
    \rlap{\kern 1.000in\lower\graphtemp\hbox to 0pt{\hss $\bullet$\hss}}%
    \special{pn 8}%
    \special{pa 1000 732}%
    \special{pa 1317 1281}%
    \special{pa 683 1281}%
    \special{pa 1000 1098}%
    \special{fp}%
    \special{pa 1317 1281}%
    \special{pa 1000 1098}%
    \special{pa 1000 732}%
    \special{pa 683 1281}%
    \special{fp}%
    \graphtemp=.5ex\advance\graphtemp by 0.659in
    \rlap{\kern 1.000in\lower\graphtemp\hbox to 0pt{\hss $y$\hss}}%
    \graphtemp=.5ex\advance\graphtemp by 1.333in
    \rlap{\kern 1.369in\lower\graphtemp\hbox to 0pt{\hss $z$\hss}}%
    \graphtemp=.5ex\advance\graphtemp by 1.333in
    \rlap{\kern 0.631in\lower\graphtemp\hbox to 0pt{\hss $x$\hss}}%
    \graphtemp=.5ex\advance\graphtemp by 1.171in
    \rlap{\kern 1.000in\lower\graphtemp\hbox to 0pt{\hss $w$\hss}}%
    \special{ar 1000 366 366 366 0 6.28319}%
    \special{ar 1634 1464 366 366 0 6.28319}%
    \special{ar 366 1464 366 366 0 6.28319}%
    \graphtemp=.5ex\advance\graphtemp by 0.366in
    \rlap{\kern 1.000in\lower\graphtemp\hbox to 0pt{\hss $Y$\hss}}%
    \graphtemp=.5ex\advance\graphtemp by 1.464in
    \rlap{\kern 1.634in\lower\graphtemp\hbox to 0pt{\hss $Z$\hss}}%
    \graphtemp=.5ex\advance\graphtemp by 1.464in
    \rlap{\kern 0.366in\lower\graphtemp\hbox to 0pt{\hss $X$\hss}}%
    \hbox{\vrule depth1.830in width0pt height 0pt}%
    \kern 2.000in
  }%
}%
}
\sp
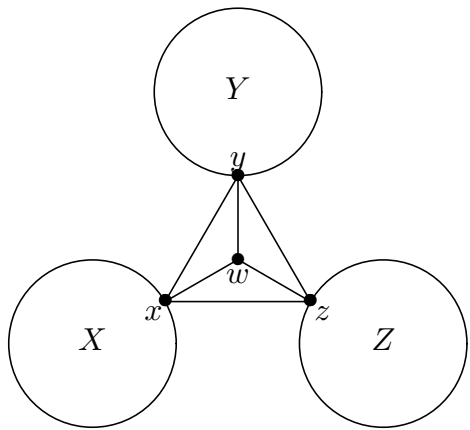
\ce{Fig.~2.  The block graph $H_k$ of diameter 3}
\sp

\TH 2.
There are trees of diameter 4 for which the bandwidth exceeds the local density
bound.
\PF
Consider the tree $T_k$ of diameter 4 illustrated in Fig. 3.  Sets $X,Y,Z$
each consist of $k-1$ leaves; $W$ consists of $k$ leaves.  Sets $X,Y,Z,W$ are
adjacent to $x,y,z,w$, respectively, and the tree is completed by making $w$
adjacent to $\{x,y,z\}$.  The tree has $4k+1$ vertices and diameter 4, with
the vertices of $X,Y,Z$ being peripheral.  The local density bound is $k$ if
$k\ge 2$, produced only by the full tree.  (When $k=1$, the tree is a star.)

Suppose that $B(H_k)=k$.  We may assume that the optimal labeling $f$ uses
labels $\{0,\ldots,4k\}$.  The distance between $f^{-1}(0)$ and $f^{-1}(4k)$
must be 4, so we may assume that $f^{-1}(0)\in X$ and $f^{-1}(4k)\in Z$.  Hence 
$f(x)=k$, $f(w)=2k$, and $f(z)=3k$.  Since $\C X=k-1$ and $\C Z=k-1$, the set
$Y\cup W$ has distinct vertices with labels less than $k$ and greater than $3k$.
Neither of these can be in $W$, since such labels differ by more than $k$ from
$f(w)$.  This yields vertices of $Y$ at distance 2 whose labels differ by more
than $2k$.
\qed

\sp
\gpic{
\expandafter\ifx\csname graph\endcsname\relax \csname newbox\endcsname\graph\fi
\expandafter\ifx\csname graphtemp\endcsname\relax \csname newdimen\endcsname\graphtemp\fi
\setbox\graph=\vtop{\vskip 0pt\hbox{%
    \graphtemp=.5ex\advance\graphtemp by 1.000in
    \rlap{\kern 0.143in\lower\graphtemp\hbox to 0pt{\hss $\bullet$\hss}}%
    \graphtemp=.5ex\advance\graphtemp by 0.429in
    \rlap{\kern 0.143in\lower\graphtemp\hbox to 0pt{\hss $\bullet$\hss}}%
    \graphtemp=.5ex\advance\graphtemp by 0.714in
    \rlap{\kern 0.429in\lower\graphtemp\hbox to 0pt{\hss $\bullet$\hss}}%
    \graphtemp=.5ex\advance\graphtemp by 0.714in
    \rlap{\kern 1.571in\lower\graphtemp\hbox to 0pt{\hss $\bullet$\hss}}%
    \graphtemp=.5ex\advance\graphtemp by 1.000in
    \rlap{\kern 1.857in\lower\graphtemp\hbox to 0pt{\hss $\bullet$\hss}}%
    \graphtemp=.5ex\advance\graphtemp by 0.429in
    \rlap{\kern 1.857in\lower\graphtemp\hbox to 0pt{\hss $\bullet$\hss}}%
    \graphtemp=.5ex\advance\graphtemp by 0.714in
    \rlap{\kern 1.000in\lower\graphtemp\hbox to 0pt{\hss $\bullet$\hss}}%
    \special{pn 8}%
    \special{pa 143 1000}%
    \special{pa 429 714}%
    \special{pa 143 429}%
    \special{pa 429 714}%
    \special{pa 1000 714}%
    \special{fp}%
    \special{pa 1857 1000}%
    \special{pa 1571 714}%
    \special{pa 1857 429}%
    \special{pa 1571 714}%
    \special{pa 1000 714}%
    \special{fp}%
    \graphtemp=.5ex\advance\graphtemp by 1.571in
    \rlap{\kern 0.714in\lower\graphtemp\hbox to 0pt{\hss $\bullet$\hss}}%
    \graphtemp=.5ex\advance\graphtemp by 1.571in
    \rlap{\kern 1.286in\lower\graphtemp\hbox to 0pt{\hss $\bullet$\hss}}%
    \graphtemp=.5ex\advance\graphtemp by 1.286in
    \rlap{\kern 1.000in\lower\graphtemp\hbox to 0pt{\hss $\bullet$\hss}}%
    \graphtemp=.5ex\advance\graphtemp by 0.143in
    \rlap{\kern 0.714in\lower\graphtemp\hbox to 0pt{\hss $\bullet$\hss}}%
    \graphtemp=.5ex\advance\graphtemp by 0.143in
    \rlap{\kern 1.000in\lower\graphtemp\hbox to 0pt{\hss $\bullet$\hss}}%
    \graphtemp=.5ex\advance\graphtemp by 0.143in
    \rlap{\kern 1.286in\lower\graphtemp\hbox to 0pt{\hss $\bullet$\hss}}%
    \special{pa 714 1571}%
    \special{pa 1000 1286}%
    \special{pa 1286 1571}%
    \special{pa 1000 1286}%
    \special{pa 1000 714}%
    \special{fp}%
    \special{pa 714 143}%
    \special{pa 1000 714}%
    \special{pa 1000 143}%
    \special{pa 1000 714}%
    \special{pa 1286 143}%
    \special{fp}%
    \special{ar 143 714 143 357 0 6.28319}%
    \special{ar 1857 714 143 357 0 6.28319}%
    \special{ar 1000 1571 357 143 0 6.28319}%
    \special{ar 1000 143 357 143 0 6.28319}%
    \graphtemp=.5ex\advance\graphtemp by 0.829in
    \rlap{\kern 0.429in\lower\graphtemp\hbox to 0pt{\hss $x$\hss}}%
    \graphtemp=.5ex\advance\graphtemp by 1.286in
    \rlap{\kern 1.114in\lower\graphtemp\hbox to 0pt{\hss $y$\hss}}%
    \graphtemp=.5ex\advance\graphtemp by 0.829in
    \rlap{\kern 1.571in\lower\graphtemp\hbox to 0pt{\hss $z$\hss}}%
    \graphtemp=.5ex\advance\graphtemp by 0.795in
    \rlap{\kern 1.081in\lower\graphtemp\hbox to 0pt{\hss $w$\hss}}%
    \graphtemp=.5ex\advance\graphtemp by 0.714in
    \rlap{\kern 0.143in\lower\graphtemp\hbox to 0pt{\hss $X$\hss}}%
    \graphtemp=.5ex\advance\graphtemp by 1.571in
    \rlap{\kern 1.000in\lower\graphtemp\hbox to 0pt{\hss $Y$\hss}}%
    \graphtemp=.5ex\advance\graphtemp by 0.714in
    \rlap{\kern 1.857in\lower\graphtemp\hbox to 0pt{\hss $Z$\hss}}%
    \graphtemp=.5ex\advance\graphtemp by 0.029in
    \rlap{\kern 1.000in\lower\graphtemp\hbox to 0pt{\hss $W$\hss}}%
    \hbox{\vrule depth1.714in width0pt height 0pt}%
    \kern 2.000in
  }%
}%
}
\sp
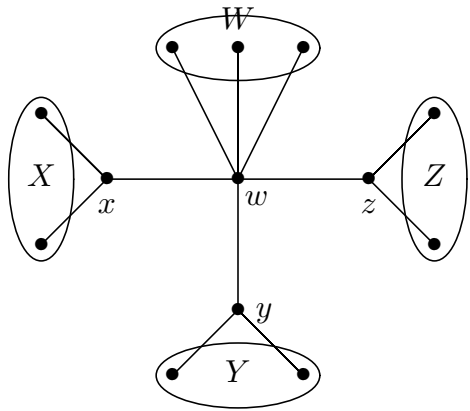
\ce{Fig.~3.  The tree $T_k$ of diameter 4}
\sp

\SH
{3. BLOCK CATERPILLARS}
We now construct optimal numberings of block caterpillars.  We view the
assignment $f$ of distinct numbers to vertices as a placement of vertices in
distinct positions; the {\it position} of $x$ is $f(x)$.  Our algorithm
constructs a numbering with minimum bandwidth, but it generally does not assign
consecutive numbers.  Condensing the vertices to consecutive positions
afterwards does not increase edge differences.  An $m$-{\it \rep\ } of a graph
(or subgraph) is a numbering such that adjacent numbers differ by at most $m$.
We use $N(S) = \UM xS N(x)$ to denote the set of vertices having a neighbor in
$S$.  A numbering $f$ is {\it faithful} if $f(x)<f(y)$ implies $f(u)<f(v)$
whenever $u,v$ are leaves adjacent to $x,y$, respectively.  We begin with two
elementary statements.

\LM 1.
If a graph $G$ has an $m$-\rep, then $G$ has a faithful $m$-\rep.
\PF
When two leaves are mis-ordered in an $m$-\rep, switching them decreases the
maximum difference on their incident edges but changes no other edge difference.
\qed

We henceforth consider only faithful numberings.  A faithful numbering of a
block graph is determined by specifying the order and position of the non-leaves
and the set of positions occupied by the leaves.

\LM 2.
Suppose $G$ is a block graph in which $X$ is a set of vertices having leaf
neighbors, and $L$ is the set of their leaf neighbors.  If $X$ occupies
consecutive positions $\al,\dots,\be$ in a faithful numbering, and the positions
of $L$ are between $\al-m$ and $\be+m$, then the differences on edges incident
to $L$ are at most $m$.
\PF
Let $X = \{x_i\}$, indexed by $f(x_i)=i$ with $\al\le i\le\be$.
Let $L = \{y_j\}$, indexed by increasing position, with $1\le j\le\C L$.
Let $u(j)=i$ if $x_i$ is the neighbor of $y_j$.  By faithfulness, $u(1)=\al$ and
$u(\C L)=\be$ and $u(j)-u(j-1)\in\{0,1\}$.  Hence $f(y_j)-f(x_{u(j)})$ is a
nondecreasing function bounded below by $-m$ and above by $m$.  \qed

\sp
For a given block caterpillar $G$, let $\VEC Q1k$ be the consecutive blocks of
the block path obtained by deleting the leaves of $G$ (in the graph of Fig. 1,
$k=4$).  Let $Q = \cup V(Q_i)$.  For $v \in Q$, let $L(v)$ denote the set of
leaves (in $G$) adjacent to $v$, and let $l(v)=\C{L(v)}$.

We first select special vertices $\SET {v_i}:{0\le i\le k+2}$.  If $k\ge 2$ and
$1< i\le k$, let $v_i$ be the shared vertex between $Q_{i-1}$ and $Q_i$.  We
may assume that $Q_1$ contains a cutvertex $v_1$ of $G$ other than $v_2$.
Otherwise, adding a leaf $x$ adjacent to a non-cutvertex of $Q_1$ yields a graph
$G'$ such that each $H'\esub G'$ containing $x$ has order and diameter one
larger than $H'-x\esub G$.  Thus $\be(G')=\be(G)$, and it suffices to study
$G'$ instead of $G$.  Similarly, we may assume that $Q_k$ contains a cutvertex
$v_{k+1}$ of $G$ other than $v_k$.  Select $v_0\in L(v_1)$ and
$v_{k+2}\in L(v_{k+1})$. 

By the same reasoning, if $k=1$ and $\C{Q_1}\ge 2$, we may assume existence of
two cutvertices $v_1,v_2$ with leaf neighbors $v_0,v_3$, respectively.  For the
degenerate case where $G$ is a star, we let $k=0$ and set $v_1$ to be the center
and $v_0,v_2$ to be arbitrary leaves.  In all cases, set $Q_0=\{v_0,v_1\}$
and $Q_{k+1}=\{v_{k+1},v_{k+2}\}$.

We will construct an $m$-\rep\ of $G$ such that $f(v_i)=im$ for $0\le i\le k+2$.
This requires putting vertices of $Q_i$ in positions $\{im,\dots,(i+1)m\}$ and
leaves adjacent to them in positions $\{(i-1)m+1,\dots,(i+2)m-1\}$
(except $\{v_0,v_{k+2}\}$).  We impose additional special properties on the
representation to facilitate the inductive argument.

\sp
\DF 1.
Let $J_i=\{im+1,\dots,(i+1)m-1\}$.  For a block caterpillar $G$ with
distinguished vertices $\VEC v0{k+2}$ as defined above, a {\it left-justified}
$m$-\rep\ is a faithful $m$-\rep\ $f$ such that
the following properties hold for $0\le i\le k+1$:
\ti 0) $f(v_i)=im$.  (Also $f(v_{k+2})=(k+2)m$.)
\ti 1) all filled positions in $J_i$ precede all unfilled positions in $J_i$.
\ti 2) if $J_i$ is not full, then $f(N(Q_i))\cap J_{i+1}=\nul$.
\ti 3) all positions for $Q_i-\{v_{i+1}\}$ precede all positions for
$L(v_{i+1})$.

We construct a left-justified $m$-\rep\ of $G$ iteratively.  The $i$th phase
produces a left-justified $m$-\rep\ of the graph $G_i$ consisting of all edges
incident to vertices of $Q_1\cup\dots\cup Q_i$.  The iteration uses the explicit
algorithm for $k=1$, so we present this as a lemma.  A block caterpillar is a
{\it clique-star} if the graph obtained by deleting all leaves is a clique; this
corresponds to $k=1$ in the description of $G$ as a block caterpillar.  When
numbering vertices, we say that an edge is {\it satisfied} if its endpoints are
at most $m$ apart.  We use $d(v)$ to denote the degree of a vertex $v$ (number
of incident edges).

\sp
\LM 3.
Every clique-star $G$ with local density at most $m$ has a left-justified
$m$-\rep.
\PF
With vertices named as above, let $X = \{\VEC x0t\}$ be the vertices of $Q$
having leaf neighbors, with $x_0=v_1$ and $x_t=v_2$, and let $Q' = Q-\{v_2\}$
and $X'=X-\{v_2\}$.  Let $f(v_i)=im$ for $0\le i\le 3$.  We will assign
positions so that $f(x_0)<\cdots<f(x_t)$ and place the leaves faithfully in
positions reserved for them.

Let $l'=\SM v{X'} l(v)$, $N=n(G)-1$, and $q=\C{Q'}$.  If $l'\le m$, reserve
positions $0,\dots,l'-1$ for leaves, and assign consecutive positions beginning
with $m$ in order to $X'$, then $Q-X$, then $L(v_2)-\{v_3\}$, skipping $2m$
(assigned to $v_2$) if $q+l(v_2)>m$.  Lemma 2 and $d(v_2)\le 2m$ imply that
every edge is satisfied, and by construction the left-justification condition
holds.

If $l'>m$, let $r$ be the least index such that $\SE j0r l(x_j)\ge m$, and let
$p=\SE j0{r-1} l(x_j)$, so $p+l(x_r)\ge m$.  Reserve all of $J_0$ for leaves.
Notice that $r<t$, by the definition of $l'$.  We consider two cases.  Each
construction fills positions other than $\{im\}$ from the left, and the
left-justification condition holds. 

If $p+l(x_r)+q\le2m$, assign positions for vertices of $Q'$ as follows: put
$\VEC x0r$ at $m,\dots,m+r$, followed by $Q-X$, and put $\VEC x{r+1}{t-1}$ at
$s-(t-r-1),\dots,s-1$, where $s=\min\{N,2m\}$.  Place the remaining leaves
in the lowest positions not yet filled.  Since $m\le p+l(X_r)\le 2m-q$, applying
Lemma 2 separately to $\VEC x0r$ and to $\VEC x{r+1}{t-1}$
guarantees that all edges are satisfied.

Finally, suppose that $p+l(x_r)+q>2m$.  This and $p<m$ force $x_r$ to have a
leaf neighbor both below $m$ and above $2m$ (see Fig.~4).  Place $\VEC x0{r-1}$
at $m,\dots,m+r-1$, and place $\VEC x{r+1}{t-1}$ at $2m-t+r+1,\dots,2m-1$.
Above $2m$, reserve the next $N-2m-1$ positions for leaves.  Reserve the
$m+t-1$ positions $m+r,\dots,2m-t+r$ for $\{x_r\}\cup(Q-X)$ and $m-q$ vertices
of $L(x_r)$.  If we can place $x_r$ in this range to satisfy the edges to its
extreme leaf neighbors, then Lemma 2 applied separately to $\VEC x0{r-1}$ and to
$\VEC x{r+1}t$ guarantees that the other edges are satisfied and completes
the proof.

To place $x_r$, observe that the lowest position in $L(x_r)$ is $p$ and the
highest is $N-p'$, where $p'=\SE i{r+1}t l(x_i)$.  Hence we need
$\max\{m+r,N-p'-m\}\le f(x_r)\le\min\{2m-t+r,p+m\}$, which requires
four inequalities.  Since $t\le q\le m$, we have $m+r\le 2m-t+r$.  Since each
vertex in $\{\VEC x0{r-1}\}$ has at least one leaf neighbor, we have
$m+r\le p+m$.  The inequality $N-p'-m\le 2m-t+r$ follows from $N\le 3m$ and
$p'\ge t-r$, which holds because each vertex in $\{\VEC x{r+1}t\}$ has at least
one leaf neighbor.  Finally, the inequality $N-p'-m\le p+m$ follows from
$N-p'-p=d(x_r)\le 2m$, which holds because $p+p'$ counts precisely the vertices
nonadjacent to $x_r$.  \qed

\sp
\gpic{
\expandafter\ifx\csname graph\endcsname\relax \csname newbox\endcsname\graph\fi
\expandafter\ifx\csname graphtemp\endcsname\relax \csname newdimen\endcsname\graphtemp\fi
\setbox\graph=\vtop{\vskip 0pt\hbox{%
    \graphtemp=.5ex\advance\graphtemp by 0.205in
    \rlap{\kern 0.058in\lower\graphtemp\hbox to 0pt{\hss $\bullet$\hss}}%
    \graphtemp=.5ex\advance\graphtemp by 0.205in
    \rlap{\kern 1.519in\lower\graphtemp\hbox to 0pt{\hss $\bullet$\hss}}%
    \graphtemp=.5ex\advance\graphtemp by 0.205in
    \rlap{\kern 2.981in\lower\graphtemp\hbox to 0pt{\hss $\bullet$\hss}}%
    \graphtemp=.5ex\advance\graphtemp by 0.205in
    \rlap{\kern 4.442in\lower\graphtemp\hbox to 0pt{\hss $\bullet$\hss}}%
    \special{pn 8}%
    \special{pa 58 205}%
    \special{pa 935 205}%
    \special{fp}%
    \special{pa 3565 205}%
    \special{pa 4222 205}%
    \special{fp}%
    \graphtemp=.5ex\advance\graphtemp by 0.058in
    \rlap{\kern 0.058in\lower\graphtemp\hbox to 0pt{\hss 0\hss}}%
    \graphtemp=.5ex\advance\graphtemp by 0.058in
    \rlap{\kern 0.877in\lower\graphtemp\hbox to 0pt{\hss $p-1$\hss}}%
    \graphtemp=.5ex\advance\graphtemp by 0.058in
    \rlap{\kern 1.519in\lower\graphtemp\hbox to 0pt{\hss $m$\hss}}%
    \graphtemp=.5ex\advance\graphtemp by 0.058in
    \rlap{\kern 2.981in\lower\graphtemp\hbox to 0pt{\hss $2m$\hss}}%
    \graphtemp=.5ex\advance\graphtemp by 0.058in
    \rlap{\kern 3.433in\lower\graphtemp\hbox to 0pt{\hss $N-p'$\hss}}%
    \graphtemp=.5ex\advance\graphtemp by 0.058in
    \rlap{\kern 4.091in\lower\graphtemp\hbox to 0pt{\hss $N-1$\hss}}%
    \graphtemp=.5ex\advance\graphtemp by 0.058in
    \rlap{\kern 4.442in\lower\graphtemp\hbox to 0pt{\hss $3m$\hss}}%
    \graphtemp=.5ex\advance\graphtemp by 0.351in
    \rlap{\kern 0.497in\lower\graphtemp\hbox to 0pt{\hss $p$ leaves\hss}}%
    \graphtemp=.5ex\advance\graphtemp by 0.351in
    \rlap{\kern 1.227in\lower\graphtemp\hbox to 0pt{\hss $L( x_r )$\hss}}%
    \graphtemp=.5ex\advance\graphtemp by 0.351in
    \rlap{\kern 1.622in\lower\graphtemp\hbox to 0pt{\hss $X'$\hss}}%
    \graphtemp=.5ex\advance\graphtemp by 0.351in
    \rlap{\kern 2.250in\lower\graphtemp\hbox to 0pt{\hss $\{x_r\}\cup (Q-X)$\hss}}%
    \graphtemp=.5ex\advance\graphtemp by 0.351in
    \rlap{\kern 2.878in\lower\graphtemp\hbox to 0pt{\hss $X'$\hss}}%
    \graphtemp=.5ex\advance\graphtemp by 0.351in
    \rlap{\kern 3.273in\lower\graphtemp\hbox to 0pt{\hss $L(x_r)$\hss}}%
    \graphtemp=.5ex\advance\graphtemp by 0.351in
    \rlap{\kern 3.879in\lower\graphtemp\hbox to 0pt{\hss $p'-1$ leaves\hss}}%
    \graphtemp=.5ex\advance\graphtemp by 0.570in
    \rlap{\kern 1.812in\lower\graphtemp\hbox to 0pt{\hss $p+l(x_r)+q+1$\hss}}%
    \special{pa 1227 570}%
    \special{pa 58 570}%
    \special{fp}%
    \special{sh 1.000}%
    \special{pa 117 584}%
    \special{pa 58 570}%
    \special{pa 117 555}%
    \special{pa 117 584}%
    \special{fp}%
    \special{pa 2396 570}%
    \special{pa 3565 570}%
    \special{fp}%
    \special{sh 1.000}%
    \special{pa 3506 555}%
    \special{pa 3565 570}%
    \special{pa 3506 584}%
    \special{pa 3506 555}%
    \special{fp}%
    \hbox{\vrule depth0.628in width0pt height 0pt}%
    \kern 4.500in
  }%
}%
}
\sp
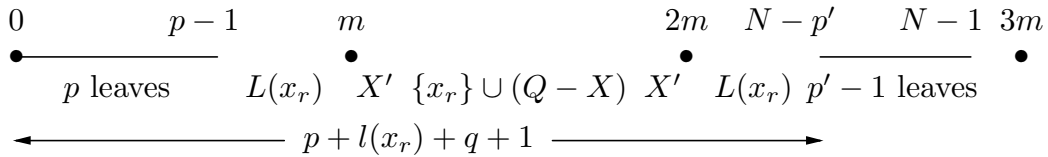
\ce{Fig.~4.  Optimal numbering of a clique-star}
\sp

\TH 3.
For every block caterpillar $G$, the bandwidth $B(G)$ equals the local density
$\be(G)$.  Furthermore, if $G$ is a block caterpillar and $\be(G)\le m$, then
$G$ has a left-justified $m$-\rep\ produced by a linear-time algorithm.
\PF
Suppose that $\diam G =k+2$.  As argued earlier, we may assume that deleting the
1-valent vertices of $G$ (except for a pair $v_0,v_{k+2}$ at maximum distance)
produces a block path with cliques $\VEC Q0{k+1}$.  Furthermore, $Q_0$
and $Q_{k+1}$ have order 2, and $\VEC v0{k+2}$ is a chordless path in $G$ having
maximum length, and $v_i$ is the cut-vertex between $Q_{i-1}$ and $Q_i$ for
$1\le i\le k+1$.

We consider a two-parameter family of subgraphs of $G$.  Let $G(h,i)$ be the
subgraph consisting of edges incident to the vertices of $Q_h\cup\dots\cup Q_i$.
For fixed $i-h\ge0$, the subgraphs $G(h,i)$ are the maximal subgraphs of
diameter $i-h+3$.  Hence if $\be_1 = \max\C{Q_i}-1$, $\be_2=\max d(v_i)/2$, and
$\be'=\max_{h\le i}\CL{\FR{n(G(h,i))-1}{i-h+3}}$, then
$\max\{\be_1,\be_2,\be'\}=\be(G)\le m$.

Let $G_i = G(1,i)$.  We produce a left-justified $m$-\rep\ of each $G_i$, by
induction on $i$, finishing with such a \rep\ for $G_k=G$.  Note that $G_{i-1}$
contains all the vertices of $Q_i-\{v_i\}$ as leaves if $2\le i\le k$.  For
$i>1$, the graph $G_i$ is obtained from $G_{i-1}$ by adding a clique on the
vertices of $Q_i-\{v_i\}$, adding the pendant edges incident to these vertices,
and adding the vertices of $Q_{i+1}-\{v_i\}$ as leaves adjacent to $v_i$.  Let
$Q_i'=Q_i-\{v_{i+1}\}$ for $1\le i\le k$.

For $i=1$, we apply Lemma 3, since $G_1$ is a clique-star.  For $i>1$, assume
that we have a left-justified $m$-\rep\ $f$ of $G_{i-1}$.  Leaf neighbors of
vertices in $Q_{i-1}'$ may locate in $J_i$ under $f$, but only if $J_{i-1}$ is
full.  Let $L'$ be the set of leaf neighbors of $Q_{i-1}'$ in positions above
$im$; these vertices do not belong to $G(i,i)$.  Nevertheless, let $G'$ be
the clique-star consisting of $G(i,i)$ together with edges from $v_i$ to $L'$
(see Fig.~5).  Because $f$ is left-justified, $L'=\nul$ if $J_{i-1}$ is not
filled by $f$.

We claim that $\be(G')\le m$.  Since $\be(G(i,i))\le m$, this fails only if
the subgraphs involving $L'$ are too big, meaning $d_{G'}(v_i)>2m$ or
$n(G')>3m+1$.  Since all neighbors of $v_i$ in $G'$ have labels between
$(i-1)m$ and $(i+1)m$ in $f$, we have $d_{G'}(v_i)\le 2m$.  To bound $n(G')$,
choose $h$ to be the largest integer in $\{1,\dots,i-1\}$ such that $J_h$ is not
full, or $h=0$ if all these intervals are full.  Since $f$ is left-justified,
every vertex of $G_{i-1}$ in a position above $hm$ belongs to $G(h+1,i)$.
Since $\VEC J{h+1}{i-1}$ are full and contain the vertices of $G(h+1,i)-G'$,
we have $n(G')=n(G(h+1,i))-(i-1-h)m$.  To bound $n(G(h+1,i))$, we use
$\diam (G(h+1,i))=i-h+2$ and $\be(G(h+1,i))\le m$ to obtain
$n(G')\le m(i-h+2)+1-(i-1-h)m= 3m+1$.

Now Lemma 3 yields a left-justified $m$-\rep\ of $G'$.  We shift each vertex
$(i-1)m$ positions rightward to obtain an $m$-\rep\ $f'$ of $G'$ using positions
between $(i-1)m$ and $(i+2)m$.  Since the vertices of
$S=Q_{i-1}'\cup L'$ are leaf neighbors of $v_i$ in $G'$, they occupy the
lowest positions under $f'$ (except that $v_i$ itself may be among them).  The
only other vertices occupying positions in both $f$ and $f'$ are those of
$T=L(v_i)\cup Q_i$, which are leaf neighbors of $v_i$ in $G_{i-1}$.  Since $f$
is left-justified, the vertices of $T$ also receive higher labels than those of
$S$ in $f$.  Hence the positions assigned to $S$ are the same in $f'$ and $f$
(those of $T$ may have moved).  We can make the vertices of $S$ occur in the
same order in $f'$ as in $f$, since these vertices are leaves in $G'$.  

We define the new $m$-\rep\ $f''$ by using $f'$ to assign positions above
$(i-1)m$ and $f$ to assign positions below $(i-1)m$.  Since $f$ and $f'$ agree
on $S$ and there are no edges from $T$ to vertices not in $G(i,i)$, we have
satisfied all edges.  The fact that $f''$ is left-justified follows from $f$ and
$f'$ being left-justified.

We comment on the complexity of the algorithm.  The graph $G$ is completely
described by giving the set of vertices in each $Q_i$ and the number of leaf
neighbors of each clique vertex.  The construction in Lemma 3 uses only these
numbers, the number of additions and subtractions involving each one is 
bounded by a constant, and the information is not used further as we proceed
in the iteration.  Thus the algorithm runs in linear time.  \qed

\sp
\gpic{
\expandafter\ifx\csname graph\endcsname\relax \csname newbox\endcsname\graph\fi
\expandafter\ifx\csname graphtemp\endcsname\relax \csname newdimen\endcsname\graphtemp\fi
\setbox\graph=\vtop{\vskip 0pt\hbox{%
    \special{pn 8}%
    \special{ar 563 328 469 234 0 6.28319}%
    \special{ar 1500 328 469 234 0 6.28319}%
    \special{ar 2438 328 469 234 0 6.28319}%
    \graphtemp=.5ex\advance\graphtemp by 0.516in
    \rlap{\kern 0.562in\lower\graphtemp\hbox to 0pt{\hss $\bullet$\hss}}%
    \graphtemp=.5ex\advance\graphtemp by 0.328in
    \rlap{\kern 0.094in\lower\graphtemp\hbox to 0pt{\hss $\bullet$\hss}}%
    \graphtemp=.5ex\advance\graphtemp by 0.141in
    \rlap{\kern 0.562in\lower\graphtemp\hbox to 0pt{\hss $\bullet$\hss}}%
    \graphtemp=.5ex\advance\graphtemp by 0.328in
    \rlap{\kern 1.031in\lower\graphtemp\hbox to 0pt{\hss $\bullet$\hss}}%
    \graphtemp=.5ex\advance\graphtemp by 0.516in
    \rlap{\kern 1.500in\lower\graphtemp\hbox to 0pt{\hss $\bullet$\hss}}%
    \graphtemp=.5ex\advance\graphtemp by 0.141in
    \rlap{\kern 1.500in\lower\graphtemp\hbox to 0pt{\hss $\bullet$\hss}}%
    \graphtemp=.5ex\advance\graphtemp by 0.328in
    \rlap{\kern 1.969in\lower\graphtemp\hbox to 0pt{\hss $\bullet$\hss}}%
    \graphtemp=.5ex\advance\graphtemp by 0.516in
    \rlap{\kern 2.438in\lower\graphtemp\hbox to 0pt{\hss $\bullet$\hss}}%
    \graphtemp=.5ex\advance\graphtemp by 0.141in
    \rlap{\kern 2.438in\lower\graphtemp\hbox to 0pt{\hss $\bullet$\hss}}%
    \special{pa 563 516}%
    \special{pa 1031 328}%
    \special{pa 94 328}%
    \special{pa 1031 328}%
    \special{pa 563 141}%
    \special{fp}%
    \special{pa 1031 328}%
    \special{pa 1500 516}%
    \special{pa 1500 141}%
    \special{pa 1969 328}%
    \special{fp}%
    \special{pa 1500 516}%
    \special{pa 1969 328}%
    \special{pa 1031 328}%
    \special{pa 1500 141}%
    \special{fp}%
    \special{pa 2438 516}%
    \special{pa 1969 328}%
    \special{pa 2438 141}%
    \special{pa 1969 328}%
    \special{pa 2906 328}%
    \special{fp}%
    \graphtemp=.5ex\advance\graphtemp by 0.000in
    \rlap{\kern 0.562in\lower\graphtemp\hbox to 0pt{\hss $Q_{i-1}$\hss}}%
    \graphtemp=.5ex\advance\graphtemp by 0.000in
    \rlap{\kern 1.500in\lower\graphtemp\hbox to 0pt{\hss $Q_i$\hss}}%
    \graphtemp=.5ex\advance\graphtemp by 0.000in
    \rlap{\kern 2.438in\lower\graphtemp\hbox to 0pt{\hss $Q_{i+1}$\hss}}%
    \graphtemp=.5ex\advance\graphtemp by 0.328in
    \rlap{\kern 2.906in\lower\graphtemp\hbox to 0pt{\hss $\bullet$\hss}}%
    \graphtemp=.5ex\advance\graphtemp by 0.797in
    \rlap{\kern 1.430in\lower\graphtemp\hbox to 0pt{\hss $\bullet$\hss}}%
    \graphtemp=.5ex\advance\graphtemp by 0.797in
    \rlap{\kern 1.570in\lower\graphtemp\hbox to 0pt{\hss $\bullet$\hss}}%
    \graphtemp=.5ex\advance\graphtemp by 0.797in
    \rlap{\kern 0.504in\lower\graphtemp\hbox to 0pt{\hss $\bullet$\hss}}%
    \graphtemp=.5ex\advance\graphtemp by 0.797in
    \rlap{\kern 0.621in\lower\graphtemp\hbox to 0pt{\hss $\bullet$\hss}}%
    \special{pa 1430 797}%
    \special{pa 1500 516}%
    \special{pa 1570 797}%
    \special{fp}%
    \special{pa 504 797}%
    \special{pa 563 516}%
    \special{pa 621 797}%
    \special{da 0.047}%
    \special{ar 563 797 176 117 0 6.28319}%
    \graphtemp=.5ex\advance\graphtemp by 0.117in
    \rlap{\kern 0.094in\lower\graphtemp\hbox to 0pt{\hss $v_{i-1}$\hss}}%
    \graphtemp=.5ex\advance\graphtemp by 0.117in
    \rlap{\kern 2.906in\lower\graphtemp\hbox to 0pt{\hss $v_{i+2}$\hss}}%
    \graphtemp=.5ex\advance\graphtemp by 0.117in
    \rlap{\kern 1.031in\lower\graphtemp\hbox to 0pt{\hss $v_i$\hss}}%
    \graphtemp=.5ex\advance\graphtemp by 0.117in
    \rlap{\kern 1.969in\lower\graphtemp\hbox to 0pt{\hss $v_{i+1}$\hss}}%
    \graphtemp=.5ex\advance\graphtemp by 0.797in
    \rlap{\kern 0.293in\lower\graphtemp\hbox to 0pt{\hss $L'$\hss}}%
    \special{pa 1031 328}%
    \special{pa 1031 563}%
    \special{pa 1031 703}%
    \special{pa 984 797}%
    \special{pa 738 797}%
    \special{sp}%
    \hbox{\vrule depth0.914in width0pt height 0pt}%
    \kern 3.000in
  }%
}%
}
\sp
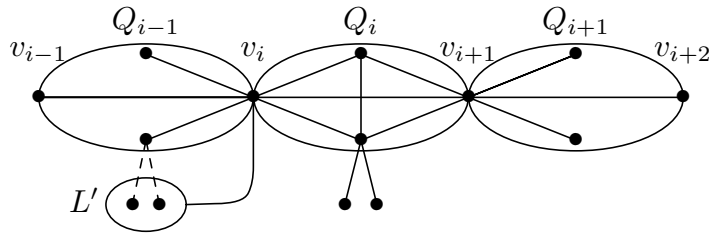
\ce{Fig.~5.  The auxiliary graph $G'$}
\sp

\SH
{4. NP-COMPLETENESS RESULTS}
Slightly enlarging the classes of 2-caterpillars or block caterpillars yields
classes on which bandwidth is NP-complete.  We prove this for two classes.
The second class consists of trees that are almost caterpillars; they have a
central path such that all other vertices except one are leaves.
The first class consists of graphs that might be called ``block
2-caterpillars'', but we use only a special subclass.  For lack of a better
name, we call the graphs in this special class ``bugs''.

\DF 2.
A graph is a {\it bug} if it is obtained from a caterpillar with an edge $xy$ on
the spine by adding a (possibly empty) clique whose vertices are adjacent to
$\{x,y\}$ and growing a nonnegative number of paths of length 2 from $x$.

\sp
All caterpillars are bugs.
Fig.~6 shows a special bug used in the NP-completeness proof.  The
{\it reflector} $R_p$ of {\it thickness} $p$ is the  bug with $5p+1$
vertices obtained from the $2p+3$-vertex caterpillar with degrees
$1,p,2,2,p,2,1$ along the spine by adding a clique of order $p-2$ adjacent to
the third edge and growing $p$ paths of length 2 from the central vertex.

\sp
\gpic{
\expandafter\ifx\csname graph\endcsname\relax \csname newbox\endcsname\graph\fi
\expandafter\ifx\csname graphtemp\endcsname\relax \csname newdimen\endcsname\graphtemp\fi
\setbox\graph=\vtop{\vskip 0pt\hbox{%
    \graphtemp=.5ex\advance\graphtemp by 0.813in
    \rlap{\kern 0.125in\lower\graphtemp\hbox to 0pt{\hss $\bullet$\hss}}%
    \graphtemp=.5ex\advance\graphtemp by 0.813in
    \rlap{\kern 0.750in\lower\graphtemp\hbox to 0pt{\hss $\bullet$\hss}}%
    \graphtemp=.5ex\advance\graphtemp by 0.813in
    \rlap{\kern 1.375in\lower\graphtemp\hbox to 0pt{\hss $\bullet$\hss}}%
    \graphtemp=.5ex\advance\graphtemp by 0.813in
    \rlap{\kern 2.000in\lower\graphtemp\hbox to 0pt{\hss $\bullet$\hss}}%
    \graphtemp=.5ex\advance\graphtemp by 0.813in
    \rlap{\kern 2.625in\lower\graphtemp\hbox to 0pt{\hss $\bullet$\hss}}%
    \graphtemp=.5ex\advance\graphtemp by 0.813in
    \rlap{\kern 3.250in\lower\graphtemp\hbox to 0pt{\hss $\bullet$\hss}}%
    \graphtemp=.5ex\advance\graphtemp by 0.813in
    \rlap{\kern 3.875in\lower\graphtemp\hbox to 0pt{\hss $\bullet$\hss}}%
    \special{pn 8}%
    \special{pa 125 813}%
    \special{pa 3875 813}%
    \special{fp}%
    \graphtemp=.5ex\advance\graphtemp by 0.688in
    \rlap{\kern 0.125in\lower\graphtemp\hbox to 0pt{\hss $a= a_0$\hss}}%
    \graphtemp=.5ex\advance\graphtemp by 0.688in
    \rlap{\kern 0.750in\lower\graphtemp\hbox to 0pt{\hss $b$\hss}}%
    \graphtemp=.5ex\advance\graphtemp by 0.688in
    \rlap{\kern 1.312in\lower\graphtemp\hbox to 0pt{\hss $c_0$\hss}}%
    \graphtemp=.5ex\advance\graphtemp by 0.688in
    \rlap{\kern 2.062in\lower\graphtemp\hbox to 0pt{\hss $w$\hss}}%
    \graphtemp=.5ex\advance\graphtemp by 0.688in
    \rlap{\kern 2.625in\lower\graphtemp\hbox to 0pt{\hss $x$\hss}}%
    \graphtemp=.5ex\advance\graphtemp by 0.688in
    \rlap{\kern 3.250in\lower\graphtemp\hbox to 0pt{\hss $y=y_0$\hss}}%
    \graphtemp=.5ex\advance\graphtemp by 0.688in
    \rlap{\kern 3.875in\lower\graphtemp\hbox to 0pt{\hss $z$\hss}}%
    \graphtemp=.5ex\advance\graphtemp by 1.438in
    \rlap{\kern 0.547in\lower\graphtemp\hbox to 0pt{\hss $\bullet$\hss}}%
    \graphtemp=.5ex\advance\graphtemp by 1.438in
    \rlap{\kern 0.953in\lower\graphtemp\hbox to 0pt{\hss $\bullet$\hss}}%
    \graphtemp=.5ex\advance\graphtemp by 1.438in
    \rlap{\kern 1.797in\lower\graphtemp\hbox to 0pt{\hss $\bullet$\hss}}%
    \graphtemp=.5ex\advance\graphtemp by 1.438in
    \rlap{\kern 2.203in\lower\graphtemp\hbox to 0pt{\hss $\bullet$\hss}}%
    \graphtemp=.5ex\advance\graphtemp by 2.062in
    \rlap{\kern 1.594in\lower\graphtemp\hbox to 0pt{\hss $\bullet$\hss}}%
    \graphtemp=.5ex\advance\graphtemp by 2.062in
    \rlap{\kern 2.406in\lower\graphtemp\hbox to 0pt{\hss $\bullet$\hss}}%
    \graphtemp=.5ex\advance\graphtemp by 1.438in
    \rlap{\kern 2.422in\lower\graphtemp\hbox to 0pt{\hss $\bullet$\hss}}%
    \graphtemp=.5ex\advance\graphtemp by 1.438in
    \rlap{\kern 2.828in\lower\graphtemp\hbox to 0pt{\hss $\bullet$\hss}}%
    \special{pa 547 1438}%
    \special{pa 750 813}%
    \special{pa 953 1438}%
    \special{fp}%
    \special{pa 1594 2063}%
    \special{pa 1797 1438}%
    \special{pa 2000 813}%
    \special{pa 2203 1438}%
    \special{pa 2406 2063}%
    \special{fp}%
    \special{pa 2422 1438}%
    \special{pa 2625 813}%
    \special{pa 2828 1438}%
    \special{fp}%
    \special{sh 0.300}%
    \special{ar 1688 344 266 219 0 6.28319}%
    \special{pa 1422 344}%
    \special{pa 1375 813}%
    \special{pa 1875 498}%
    \special{fp}%
    \special{pa 1953 344}%
    \special{pa 2000 813}%
    \special{pa 1500 498}%
    \special{fp}%
    \graphtemp=.5ex\advance\graphtemp by 1.594in
    \rlap{\kern 0.781in\lower\graphtemp\hbox to 0pt{\hss $a_1\cdots a_{p-2}$\hss}}%
    \graphtemp=.5ex\advance\graphtemp by 1.594in
    \rlap{\kern 2.687in\lower\graphtemp\hbox to 0pt{\hss $y_1\cdots y_{p-2}$\hss}}%
    \graphtemp=.5ex\advance\graphtemp by 1.656in
    \rlap{\kern 2.000in\lower\graphtemp\hbox to 0pt{\hss $w_1 \cdots w_p$\hss}}%
    \graphtemp=.5ex\advance\graphtemp by 2.062in
    \rlap{\kern 2.000in\lower\graphtemp\hbox to 0pt{\hss $w_1'\cdots w_p'$\hss}}%
    \graphtemp=.5ex\advance\graphtemp by 0.000in
    \rlap{\kern 1.719in\lower\graphtemp\hbox to 0pt{\hss $c_1\cdots c_{p-2}$\hss}}%
    \hbox{\vrule depth2.188in width0pt height 0pt}%
    \kern 4.000in
  }%
}%
}
\sp
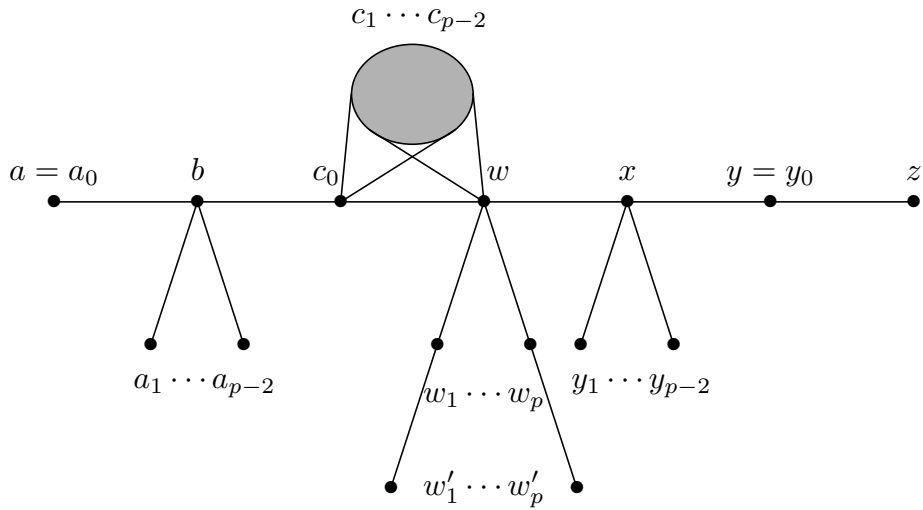
\ce{Fig.~6.  The reflector $R_p$ of thickness $p$}
\sp

We follow the method used by Monien [12] to prove that bandwidth is NP-complete
for 3-caterpillars.  A {\it peripheral vertex} of a graph is a vertex of 
maximum eccentricity, where the {\it eccentricity} of a vertex is its maximum
distance from other vertices.  Deleting the peripheral vertices of $R_p$ yields
a bug of diameter 4 with $4p+1$ vertices, so $R_p$ has local density at least
$p$, and the numbering in Fig.~7 shows that its bandwidth is $p$.

\sp
\gpic{
\expandafter\ifx\csname graph\endcsname\relax \csname newbox\endcsname\graph\fi
\expandafter\ifx\csname graphtemp\endcsname\relax \csname newdimen\endcsname\graphtemp\fi
\setbox\graph=\vtop{\vskip 0pt\hbox{%
    \graphtemp=.5ex\advance\graphtemp by 0.117in
    \rlap{\kern 0.083in\lower\graphtemp\hbox to 0pt{\hss $\bullet$\hss}}%
    \graphtemp=.5ex\advance\graphtemp by 0.117in
    \rlap{\kern 0.289in\lower\graphtemp\hbox to 0pt{\hss $\bullet$\hss}}%
    \graphtemp=.5ex\advance\graphtemp by 0.117in
    \rlap{\kern 0.495in\lower\graphtemp\hbox to 0pt{\hss $\bullet$\hss}}%
    \graphtemp=.5ex\advance\graphtemp by 0.117in
    \rlap{\kern 0.908in\lower\graphtemp\hbox to 0pt{\hss $\bullet$\hss}}%
    \graphtemp=.5ex\advance\graphtemp by 0.117in
    \rlap{\kern 1.115in\lower\graphtemp\hbox to 0pt{\hss $\bullet$\hss}}%
    \graphtemp=.5ex\advance\graphtemp by 0.117in
    \rlap{\kern 1.321in\lower\graphtemp\hbox to 0pt{\hss $\bullet$\hss}}%
    \graphtemp=.5ex\advance\graphtemp by 0.117in
    \rlap{\kern 1.528in\lower\graphtemp\hbox to 0pt{\hss $\bullet$\hss}}%
    \graphtemp=.5ex\advance\graphtemp by 0.117in
    \rlap{\kern 1.940in\lower\graphtemp\hbox to 0pt{\hss $\bullet$\hss}}%
    \graphtemp=.5ex\advance\graphtemp by 0.117in
    \rlap{\kern 2.147in\lower\graphtemp\hbox to 0pt{\hss $\bullet$\hss}}%
    \graphtemp=.5ex\advance\graphtemp by 0.241in
    \rlap{\kern 0.083in\lower\graphtemp\hbox to 0pt{\hss $z$\hss}}%
    \graphtemp=.5ex\advance\graphtemp by 0.241in
    \rlap{\kern 0.289in\lower\graphtemp\hbox to 0pt{\hss $a$\hss}}%
    \graphtemp=.5ex\advance\graphtemp by 0.241in
    \rlap{\kern 0.495in\lower\graphtemp\hbox to 0pt{\hss $a_1$\hss}}%
    \graphtemp=.5ex\advance\graphtemp by 0.138in
    \rlap{\kern 0.702in\lower\graphtemp\hbox to 0pt{\hss $\cdots$\hss}}%
    \graphtemp=.5ex\advance\graphtemp by 0.241in
    \rlap{\kern 0.908in\lower\graphtemp\hbox to 0pt{\hss $a_{p-2}$\hss}}%
    \graphtemp=.5ex\advance\graphtemp by 0.241in
    \rlap{\kern 1.115in\lower\graphtemp\hbox to 0pt{\hss $y$\hss}}%
    \graphtemp=.5ex\advance\graphtemp by 0.241in
    \rlap{\kern 1.321in\lower\graphtemp\hbox to 0pt{\hss $b$\hss}}%
    \graphtemp=.5ex\advance\graphtemp by 0.241in
    \rlap{\kern 1.528in\lower\graphtemp\hbox to 0pt{\hss $y_1$\hss}}%
    \graphtemp=.5ex\advance\graphtemp by 0.138in
    \rlap{\kern 1.734in\lower\graphtemp\hbox to 0pt{\hss $\cdots$\hss}}%
    \graphtemp=.5ex\advance\graphtemp by 0.241in
    \rlap{\kern 1.940in\lower\graphtemp\hbox to 0pt{\hss $y_{p-2}$\hss}}%
    \graphtemp=.5ex\advance\graphtemp by 0.241in
    \rlap{\kern 2.147in\lower\graphtemp\hbox to 0pt{\hss $x$\hss}}%
    \graphtemp=.5ex\advance\graphtemp by 0.241in
    \rlap{\kern 2.353in\lower\graphtemp\hbox to 0pt{\hss $c_0$\hss}}%
    \graphtemp=.5ex\advance\graphtemp by 0.241in
    \rlap{\kern 2.560in\lower\graphtemp\hbox to 0pt{\hss $c_1$\hss}}%
    \graphtemp=.5ex\advance\graphtemp by 0.117in
    \rlap{\kern 2.353in\lower\graphtemp\hbox to 0pt{\hss $\bullet$\hss}}%
    \graphtemp=.5ex\advance\graphtemp by 0.117in
    \rlap{\kern 2.560in\lower\graphtemp\hbox to 0pt{\hss $\bullet$\hss}}%
    \graphtemp=.5ex\advance\graphtemp by 0.117in
    \rlap{\kern 2.972in\lower\graphtemp\hbox to 0pt{\hss $\bullet$\hss}}%
    \graphtemp=.5ex\advance\graphtemp by 0.117in
    \rlap{\kern 3.179in\lower\graphtemp\hbox to 0pt{\hss $\bullet$\hss}}%
    \graphtemp=.5ex\advance\graphtemp by 0.117in
    \rlap{\kern 3.385in\lower\graphtemp\hbox to 0pt{\hss $\bullet$\hss}}%
    \graphtemp=.5ex\advance\graphtemp by 0.117in
    \rlap{\kern 3.798in\lower\graphtemp\hbox to 0pt{\hss $\bullet$\hss}}%
    \graphtemp=.5ex\advance\graphtemp by 0.117in
    \rlap{\kern 4.005in\lower\graphtemp\hbox to 0pt{\hss $\bullet$\hss}}%
    \graphtemp=.5ex\advance\graphtemp by 0.117in
    \rlap{\kern 4.417in\lower\graphtemp\hbox to 0pt{\hss $\bullet$\hss}}%
    \graphtemp=.5ex\advance\graphtemp by 0.138in
    \rlap{\kern 2.766in\lower\graphtemp\hbox to 0pt{\hss $\cdots$\hss}}%
    \graphtemp=.5ex\advance\graphtemp by 0.241in
    \rlap{\kern 2.972in\lower\graphtemp\hbox to 0pt{\hss $c_{p-2}$\hss}}%
    \graphtemp=.5ex\advance\graphtemp by 0.241in
    \rlap{\kern 3.179in\lower\graphtemp\hbox to 0pt{\hss $w$\hss}}%
    \graphtemp=.5ex\advance\graphtemp by 0.241in
    \rlap{\kern 3.385in\lower\graphtemp\hbox to 0pt{\hss $w_1$\hss}}%
    \graphtemp=.5ex\advance\graphtemp by 0.138in
    \rlap{\kern 3.592in\lower\graphtemp\hbox to 0pt{\hss $\cdots$\hss}}%
    \graphtemp=.5ex\advance\graphtemp by 0.241in
    \rlap{\kern 3.798in\lower\graphtemp\hbox to 0pt{\hss $w_p$\hss}}%
    \graphtemp=.5ex\advance\graphtemp by 0.241in
    \rlap{\kern 4.005in\lower\graphtemp\hbox to 0pt{\hss $w_1'$\hss}}%
    \graphtemp=.5ex\advance\graphtemp by 0.138in
    \rlap{\kern 4.211in\lower\graphtemp\hbox to 0pt{\hss $\cdots$\hss}}%
    \graphtemp=.5ex\advance\graphtemp by 0.241in
    \rlap{\kern 4.417in\lower\graphtemp\hbox to 0pt{\hss $w_p'$\hss}}%
    \special{pn 8}%
    \special{ar 599 1197 1197 1197 -2.016435 -1.125157}%
    \special{ar 805 1197 1197 1197 -2.016435 -1.125157}%
    \special{ar 908 1040 1011 1011 -1.991238 -1.150355}%
    \special{ar 1115 679 599 599 -1.922851 -1.218741}%
    \special{ar 1631 1197 1197 1197 -2.016435 -1.125157}%
    \special{ar 1837 860 805 805 -1.965587 -1.176005}%
    \special{ar 2044 474 372 372 -1.852276 -1.289316}%
    \special{ar 1837 1197 1197 1197 -2.016435 -1.125157}%
    \special{ar 2663 1197 1197 1197 -2.016435 -1.125157}%
    \special{ar 2766 1040 1011 1011 -1.991238 -1.150355}%
    \special{ar 2869 860 805 805 -1.965587 -1.176005}%
    \special{ar 3076 474 372 372 -1.852276 -1.289316}%
    \special{ar 3282 474 372 372 -1.852276 -1.289316}%
    \special{ar 3489 860 805 805 -1.965587 -1.176005}%
    \special{ar 3695 860 805 805 -1.965587 -1.176005}%
    \special{ar 4108 860 805 805 -1.965587 -1.176005}%
    \hbox{\vrule depth0.406in width0pt height 0pt}%
    \kern 4.500in
  }%
}%
}
\sp
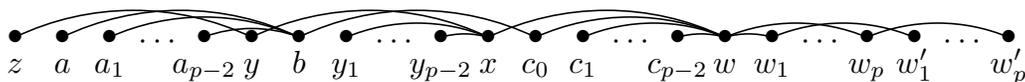
\ce{Fig.~7.  An optimal numbering of $R_p$}
\sp

The property of $R_p$ needed for the NP-completeness reduction is that in every
optimal numbering of $R_p$, the peripheral vertices appear on the same end.  As
a result, placing the reflector in the middle of a caterpillar-like object
forces its vertices to take positions on one end of a numbering that achieves a
desired bandwidth.  This motivates the term ``reflector''.

\LM 4.
Let $R_p$ be the reflector of thickness $p$, where $p\ge4$.
In every optimal numbering of $R_p$ with positions $\{0,\dots,5p\}$,
the positions of the peripheral vertices $a$ and $z$ are both below $p$ or
both above $4p$.
\PF
Let $f$ be an optimal numbering of $R_p$, so $B(f)=p$.  By symmetry, we may
assume that $f(x)<f(w)$.  Since $R_p$ has $4p+1$ vertices with distance at most
2 from $w$, there are at least $2p$ positions to each side of $f(w)$; thus
$2p\le f(w)\le3p]$.  Since $d(w)=2p$, vertices at distance 2 from $w$ occupy
positions more than $p$ from $f(w)$.  Since $f(x)<f(w)$, this forces
$f(y_0),\dots,f(y_{p-2})$ into the interval $[f(w)-2p,f(w)-p-1]$.

If this interval also contains $f(b)$, then $f(a)$ and $f(z)$ must both be below
$f(w)-2p$, since no other positions remain within $p$ of $f(b)$ or $f(y)$.
Since $f(w)\le 3p$, this will complete the proof.

Since $d(b,w)=2$, it now suffices to show that $f(b)<f(w)$.  For each vertex $v$
at distance 2 from $w$, $p<|f(v)-f(w)|\le2p$.  With $f(y_0),\dots,f(y_{p-2})$
below $f(w)$, this forces at least $p-1$ of $f(w_1'),\dots,f(w_p')$ into
$[f(w)+p+1,f(w)+2p]$.  This in turn forces at least $p-1$ of
$f(w_1),\dots,f(w_p)$ into $[f(w)+1,f(w)+p]$.  At most one position remains in
this interval, so at least $p-2$ of $f(c_0),\dots,f(c_{p-2})$ lie in
$[f(w)-p,f(w)-1]$.  All of these $p-2$ vertices are within distance 2 from $b$
in $R_p$, and thus $f(b)\le(f(w)-p+2)+2p$.

If $f(b)>f(w)$, then $f(b)\in\{f(w)+p+1,f(w)+p+2\}$,
because the neighbors of $w$ occupy the positions within $p$ of $f(w)$.
This forces $f(a_0),\dots,f(a_{p-2})$ to occupy positions above
$2p$.  When $p\ge4$, there are at least 3 such values, and one of them
is now too far above $f(b)$.

When $f(x)>f(w)$, the analogous argument locates $a$ and $z$ above $4p$.
\qed

We prove NP-completeness of the bandwidth problem for bugs by reduction from the
Multiprocessor Scheduling Problem.  An instance $T$ of this problem consists of
a number $m$ of processors, a deadline $D$, and $n$ tasks with integer execution
times $\VEC t1n$.  The decision problem asks whether the tasks can be assigned
to the processors such that for each processor, the total execution time for the
assigned tasks is at most $D$.  When the answer is ``Yes'', we say that the
instance is {\it solvable}.  As shown in [5, p95-106], Multiprocessor
Scheduling is NP-complete in the strong sense, which means that we can consider
the size of $T$ to be $m+n+\max t_i$.

\TH 4.
The bandwidth problem is NP-complete for bugs.
\PF
It suffices to show that for each instance $T=(m,D,\VEC t1n)$ of the
Multiprocessor Scheduling Problem, we can construct a bug $G$ and an integer
$b$ such that $T$ is solvable if and only if $B(G)\le b$.
Furthermore, the construction must run in time polynomial in $m$, $n$, and
$\max t_i$.

Given $T$, choose $p$ such that $p>2n(D+4)$.  Let $b=p+1+2n$, and let
$D'=2m(D+2)-4$.  We construct a bug $G$ using two caterpillars and the reflector
$R_b$ of thickness $b$ (see Fig.~8).  Caterpillar $C$ consists of
$m(D+2+2p)+4n$ vertices, with $\lmb=m(D+2)$ vertices on the spine.  The
$2mp+4n$ additional leaves of $C$ are attached as follows: $2p+4n$ at the
second vertex and $2p$ at each vertex whose distance from the second vertex
along the spine is a multiple of $D+2$.  Caterpillar $C'$ consists of
$(p\SE i1n t_i)+nD'$ vertices, with $\SE i1n t_i +nD'$ vertices on the spine.
The $i$th {\it segment} of $C'$ consists of $t_i$ vertices each with $p-1$
leaves as neighbors followed by $D'$ vertices with no leaf neighbors; the 
$t_i p$ vertices appearing first are the {\it task vertices}.  To complete the
bug $G$, add two edges: one each from the peripheral vertices $a$ and $z$ of
$R_p$ to the last vertices on the spines of $C$ and $C'$, respectively.
Note that the number of vertices in $G$ is given by a polynomial in $n$, $D$,
and $\SE i1n t_i$.

\sp
\gpic{
\expandafter\ifx\csname graph\endcsname\relax \csname newbox\endcsname\graph\fi
\expandafter\ifx\csname graphtemp\endcsname\relax \csname newdimen\endcsname\graphtemp\fi
\setbox\graph=\vtop{\vskip 0pt\hbox{%
    \graphtemp=.5ex\advance\graphtemp by 1.540in
    \rlap{\kern 0.183in\lower\graphtemp\hbox to 0pt{\hss $\bullet$\hss}}%
    \graphtemp=.5ex\advance\graphtemp by 1.540in
    \rlap{\kern 0.391in\lower\graphtemp\hbox to 0pt{\hss $\bullet$\hss}}%
    \graphtemp=.5ex\advance\graphtemp by 1.540in
    \rlap{\kern 0.599in\lower\graphtemp\hbox to 0pt{\hss $\bullet$\hss}}%
    \graphtemp=.5ex\advance\graphtemp by 1.540in
    \rlap{\kern 0.807in\lower\graphtemp\hbox to 0pt{\hss $\bullet$\hss}}%
    \graphtemp=.5ex\advance\graphtemp by 1.540in
    \rlap{\kern 1.223in\lower\graphtemp\hbox to 0pt{\hss $\bullet$\hss}}%
    \graphtemp=.5ex\advance\graphtemp by 1.956in
    \rlap{\kern 0.260in\lower\graphtemp\hbox to 0pt{\hss $\bullet$\hss}}%
    \graphtemp=.5ex\advance\graphtemp by 1.956in
    \rlap{\kern 0.522in\lower\graphtemp\hbox to 0pt{\hss $\bullet$\hss}}%
    \special{pn 8}%
    \special{pa 183 1540}%
    \special{pa 870 1540}%
    \special{fp}%
    \special{pa 1161 1540}%
    \special{pa 1348 1540}%
    \special{fp}%
    \special{pa 260 1956}%
    \special{pa 391 1540}%
    \special{fp}%
    \special{pa 391 1540}%
    \special{pa 522 1956}%
    \special{fp}%
    \graphtemp=.5ex\advance\graphtemp by 1.560in
    \rlap{\kern 1.015in\lower\graphtemp\hbox to 0pt{\hss $\cdots$\hss}}%
    \special{pa 807 1644}%
    \special{pa 807 1748}%
    \special{pa 849 1748}%
    \special{fp}%
    \special{pa 1223 1644}%
    \special{pa 1223 1748}%
    \special{pa 1161 1748}%
    \special{fp}%
    \graphtemp=.5ex\advance\graphtemp by 1.976in
    \rlap{\kern 0.391in\lower\graphtemp\hbox to 0pt{\hss $\cdots$\hss}}%
    \graphtemp=.5ex\advance\graphtemp by 2.080in
    \rlap{\kern 0.391in\lower\graphtemp\hbox to 0pt{\hss $2p+4n$\hss}}%
    \graphtemp=.5ex\advance\graphtemp by 1.748in
    \rlap{\kern 1.015in\lower\graphtemp\hbox to 0pt{\hss $D-1$\hss}}%
    \graphtemp=.5ex\advance\graphtemp by 1.540in
    \rlap{\kern 1.431in\lower\graphtemp\hbox to 0pt{\hss $\bullet$\hss}}%
    \graphtemp=.5ex\advance\graphtemp by 1.540in
    \rlap{\kern 1.639in\lower\graphtemp\hbox to 0pt{\hss $\bullet$\hss}}%
    \graphtemp=.5ex\advance\graphtemp by 1.540in
    \rlap{\kern 1.847in\lower\graphtemp\hbox to 0pt{\hss $\bullet$\hss}}%
    \graphtemp=.5ex\advance\graphtemp by 1.540in
    \rlap{\kern 2.055in\lower\graphtemp\hbox to 0pt{\hss $\bullet$\hss}}%
    \graphtemp=.5ex\advance\graphtemp by 1.540in
    \rlap{\kern 2.472in\lower\graphtemp\hbox to 0pt{\hss $\bullet$\hss}}%
    \graphtemp=.5ex\advance\graphtemp by 1.956in
    \rlap{\kern 1.508in\lower\graphtemp\hbox to 0pt{\hss $\bullet$\hss}}%
    \graphtemp=.5ex\advance\graphtemp by 1.956in
    \rlap{\kern 1.770in\lower\graphtemp\hbox to 0pt{\hss $\bullet$\hss}}%
    \special{pa 1431 1540}%
    \special{pa 2118 1540}%
    \special{fp}%
    \special{pa 2409 1540}%
    \special{pa 2596 1540}%
    \special{fp}%
    \special{pa 1508 1956}%
    \special{pa 1639 1540}%
    \special{fp}%
    \special{pa 1639 1540}%
    \special{pa 1770 1956}%
    \special{fp}%
    \graphtemp=.5ex\advance\graphtemp by 1.560in
    \rlap{\kern 2.264in\lower\graphtemp\hbox to 0pt{\hss $\cdots$\hss}}%
    \special{pa 2055 1644}%
    \special{pa 2055 1748}%
    \special{pa 2097 1748}%
    \special{fp}%
    \special{pa 2472 1644}%
    \special{pa 2472 1748}%
    \special{pa 2409 1748}%
    \special{fp}%
    \graphtemp=.5ex\advance\graphtemp by 1.976in
    \rlap{\kern 1.639in\lower\graphtemp\hbox to 0pt{\hss $\cdots$\hss}}%
    \graphtemp=.5ex\advance\graphtemp by 2.080in
    \rlap{\kern 1.639in\lower\graphtemp\hbox to 0pt{\hss $2p$\hss}}%
    \graphtemp=.5ex\advance\graphtemp by 1.748in
    \rlap{\kern 2.264in\lower\graphtemp\hbox to 0pt{\hss $D-1$\hss}}%
    \graphtemp=.5ex\advance\graphtemp by 1.540in
    \rlap{\kern 3.096in\lower\graphtemp\hbox to 0pt{\hss $\bullet$\hss}}%
    \graphtemp=.5ex\advance\graphtemp by 1.540in
    \rlap{\kern 3.304in\lower\graphtemp\hbox to 0pt{\hss $\bullet$\hss}}%
    \graphtemp=.5ex\advance\graphtemp by 1.540in
    \rlap{\kern 3.512in\lower\graphtemp\hbox to 0pt{\hss $\bullet$\hss}}%
    \graphtemp=.5ex\advance\graphtemp by 1.540in
    \rlap{\kern 3.720in\lower\graphtemp\hbox to 0pt{\hss $\bullet$\hss}}%
    \graphtemp=.5ex\advance\graphtemp by 1.540in
    \rlap{\kern 4.136in\lower\graphtemp\hbox to 0pt{\hss $\bullet$\hss}}%
    \graphtemp=.5ex\advance\graphtemp by 1.956in
    \rlap{\kern 3.173in\lower\graphtemp\hbox to 0pt{\hss $\bullet$\hss}}%
    \graphtemp=.5ex\advance\graphtemp by 1.956in
    \rlap{\kern 3.435in\lower\graphtemp\hbox to 0pt{\hss $\bullet$\hss}}%
    \special{pa 3096 1540}%
    \special{pa 3782 1540}%
    \special{fp}%
    \special{pa 4074 1540}%
    \special{pa 4136 1540}%
    \special{fp}%
    \special{pa 3173 1956}%
    \special{pa 3304 1540}%
    \special{fp}%
    \special{pa 3304 1540}%
    \special{pa 3435 1956}%
    \special{fp}%
    \graphtemp=.5ex\advance\graphtemp by 1.560in
    \rlap{\kern 3.928in\lower\graphtemp\hbox to 0pt{\hss $\cdots$\hss}}%
    \special{pa 3720 1644}%
    \special{pa 3720 1748}%
    \special{pa 3761 1748}%
    \special{fp}%
    \special{pa 4136 1644}%
    \special{pa 4136 1748}%
    \special{pa 4074 1748}%
    \special{fp}%
    \graphtemp=.5ex\advance\graphtemp by 1.976in
    \rlap{\kern 3.304in\lower\graphtemp\hbox to 0pt{\hss $\cdots$\hss}}%
    \graphtemp=.5ex\advance\graphtemp by 2.080in
    \rlap{\kern 3.304in\lower\graphtemp\hbox to 0pt{\hss $2p$\hss}}%
    \graphtemp=.5ex\advance\graphtemp by 1.748in
    \rlap{\kern 3.928in\lower\graphtemp\hbox to 0pt{\hss $D-1$\hss}}%
    \graphtemp=.5ex\advance\graphtemp by 1.560in
    \rlap{\kern 2.888in\lower\graphtemp\hbox to 0pt{\hss $\cdots$\hss}}%
    \graphtemp=.5ex\advance\graphtemp by 1.956in
    \rlap{\kern 2.472in\lower\graphtemp\hbox to 0pt{\hss $m$ segments\hss}}%
    \graphtemp=.5ex\advance\graphtemp by 2.164in
    \rlap{\kern 2.472in\lower\graphtemp\hbox to 0pt{\hss $C$\hss}}%
    \graphtemp=.5ex\advance\graphtemp by 0.291in
    \rlap{\kern 0.183in\lower\graphtemp\hbox to 0pt{\hss $\bullet$\hss}}%
    \graphtemp=.5ex\advance\graphtemp by 0.291in
    \rlap{\kern 0.599in\lower\graphtemp\hbox to 0pt{\hss $\bullet$\hss}}%
    \graphtemp=.5ex\advance\graphtemp by 0.291in
    \rlap{\kern 0.807in\lower\graphtemp\hbox to 0pt{\hss $\bullet$\hss}}%
    \graphtemp=.5ex\advance\graphtemp by 0.291in
    \rlap{\kern 1.223in\lower\graphtemp\hbox to 0pt{\hss $\bullet$\hss}}%
    \graphtemp=.5ex\advance\graphtemp by 0.707in
    \rlap{\kern 0.052in\lower\graphtemp\hbox to 0pt{\hss $\bullet$\hss}}%
    \graphtemp=.5ex\advance\graphtemp by 0.707in
    \rlap{\kern 0.314in\lower\graphtemp\hbox to 0pt{\hss $\bullet$\hss}}%
    \graphtemp=.5ex\advance\graphtemp by 0.707in
    \rlap{\kern 0.468in\lower\graphtemp\hbox to 0pt{\hss $\bullet$\hss}}%
    \graphtemp=.5ex\advance\graphtemp by 0.707in
    \rlap{\kern 0.730in\lower\graphtemp\hbox to 0pt{\hss $\bullet$\hss}}%
    \special{pa 52 707}%
    \special{pa 183 291}%
    \special{pa 314 707}%
    \special{fp}%
    \special{pa 468 707}%
    \special{pa 599 291}%
    \special{pa 730 707}%
    \special{fp}%
    \special{pa 183 291}%
    \special{pa 260 291}%
    \special{fp}%
    \special{pa 522 291}%
    \special{pa 884 291}%
    \special{fp}%
    \special{pa 1146 291}%
    \special{pa 1348 291}%
    \special{fp}%
    \graphtemp=.5ex\advance\graphtemp by 0.728in
    \rlap{\kern 0.183in\lower\graphtemp\hbox to 0pt{\hss $\cdots$\hss}}%
    \graphtemp=.5ex\advance\graphtemp by 0.728in
    \rlap{\kern 0.599in\lower\graphtemp\hbox to 0pt{\hss $\cdots$\hss}}%
    \graphtemp=.5ex\advance\graphtemp by 0.832in
    \rlap{\kern 0.183in\lower\graphtemp\hbox to 0pt{\hss $p-1$\hss}}%
    \graphtemp=.5ex\advance\graphtemp by 0.832in
    \rlap{\kern 0.599in\lower\graphtemp\hbox to 0pt{\hss $p-1$\hss}}%
    \graphtemp=.5ex\advance\graphtemp by 0.083in
    \rlap{\kern 0.391in\lower\graphtemp\hbox to 0pt{\hss $t_1$\hss}}%
    \graphtemp=.5ex\advance\graphtemp by 0.104in
    \rlap{\kern 1.015in\lower\graphtemp\hbox to 0pt{\hss $D'$\hss}}%
    \graphtemp=.5ex\advance\graphtemp by 0.312in
    \rlap{\kern 0.391in\lower\graphtemp\hbox to 0pt{\hss $\cdots$\hss}}%
    \graphtemp=.5ex\advance\graphtemp by 0.312in
    \rlap{\kern 1.015in\lower\graphtemp\hbox to 0pt{\hss $\cdots$\hss}}%
    \special{pa 183 187}%
    \special{pa 183 83}%
    \special{pa 225 83}%
    \special{fp}%
    \special{pa 599 187}%
    \special{pa 599 83}%
    \special{pa 537 83}%
    \special{fp}%
    \special{pa 807 187}%
    \special{pa 807 83}%
    \special{pa 849 83}%
    \special{fp}%
    \special{pa 1223 187}%
    \special{pa 1223 83}%
    \special{pa 1161 83}%
    \special{fp}%
    \graphtemp=.5ex\advance\graphtemp by 0.291in
    \rlap{\kern 1.431in\lower\graphtemp\hbox to 0pt{\hss $\bullet$\hss}}%
    \graphtemp=.5ex\advance\graphtemp by 0.291in
    \rlap{\kern 1.847in\lower\graphtemp\hbox to 0pt{\hss $\bullet$\hss}}%
    \graphtemp=.5ex\advance\graphtemp by 0.291in
    \rlap{\kern 2.055in\lower\graphtemp\hbox to 0pt{\hss $\bullet$\hss}}%
    \graphtemp=.5ex\advance\graphtemp by 0.291in
    \rlap{\kern 2.472in\lower\graphtemp\hbox to 0pt{\hss $\bullet$\hss}}%
    \graphtemp=.5ex\advance\graphtemp by 0.707in
    \rlap{\kern 1.300in\lower\graphtemp\hbox to 0pt{\hss $\bullet$\hss}}%
    \graphtemp=.5ex\advance\graphtemp by 0.707in
    \rlap{\kern 1.562in\lower\graphtemp\hbox to 0pt{\hss $\bullet$\hss}}%
    \graphtemp=.5ex\advance\graphtemp by 0.707in
    \rlap{\kern 1.716in\lower\graphtemp\hbox to 0pt{\hss $\bullet$\hss}}%
    \graphtemp=.5ex\advance\graphtemp by 0.707in
    \rlap{\kern 1.979in\lower\graphtemp\hbox to 0pt{\hss $\bullet$\hss}}%
    \special{pa 1300 707}%
    \special{pa 1431 291}%
    \special{pa 1562 707}%
    \special{fp}%
    \special{pa 1716 707}%
    \special{pa 1847 291}%
    \special{pa 1979 707}%
    \special{fp}%
    \special{pa 1431 291}%
    \special{pa 1508 291}%
    \special{fp}%
    \special{pa 1770 291}%
    \special{pa 2132 291}%
    \special{fp}%
    \special{pa 2395 291}%
    \special{pa 2596 291}%
    \special{fp}%
    \graphtemp=.5ex\advance\graphtemp by 0.728in
    \rlap{\kern 1.431in\lower\graphtemp\hbox to 0pt{\hss $\cdots$\hss}}%
    \graphtemp=.5ex\advance\graphtemp by 0.728in
    \rlap{\kern 1.847in\lower\graphtemp\hbox to 0pt{\hss $\cdots$\hss}}%
    \graphtemp=.5ex\advance\graphtemp by 0.832in
    \rlap{\kern 1.431in\lower\graphtemp\hbox to 0pt{\hss $p-1$\hss}}%
    \graphtemp=.5ex\advance\graphtemp by 0.832in
    \rlap{\kern 1.847in\lower\graphtemp\hbox to 0pt{\hss $p-1$\hss}}%
    \graphtemp=.5ex\advance\graphtemp by 0.083in
    \rlap{\kern 1.639in\lower\graphtemp\hbox to 0pt{\hss $t_2$\hss}}%
    \graphtemp=.5ex\advance\graphtemp by 0.104in
    \rlap{\kern 2.264in\lower\graphtemp\hbox to 0pt{\hss $D'$\hss}}%
    \graphtemp=.5ex\advance\graphtemp by 0.312in
    \rlap{\kern 1.639in\lower\graphtemp\hbox to 0pt{\hss $\cdots$\hss}}%
    \graphtemp=.5ex\advance\graphtemp by 0.312in
    \rlap{\kern 2.264in\lower\graphtemp\hbox to 0pt{\hss $\cdots$\hss}}%
    \special{pa 1431 187}%
    \special{pa 1431 83}%
    \special{pa 1473 83}%
    \special{fp}%
    \special{pa 1847 187}%
    \special{pa 1847 83}%
    \special{pa 1785 83}%
    \special{fp}%
    \special{pa 2055 187}%
    \special{pa 2055 83}%
    \special{pa 2097 83}%
    \special{fp}%
    \special{pa 2472 187}%
    \special{pa 2472 83}%
    \special{pa 2409 83}%
    \special{fp}%
    \graphtemp=.5ex\advance\graphtemp by 0.291in
    \rlap{\kern 3.096in\lower\graphtemp\hbox to 0pt{\hss $\bullet$\hss}}%
    \graphtemp=.5ex\advance\graphtemp by 0.291in
    \rlap{\kern 3.512in\lower\graphtemp\hbox to 0pt{\hss $\bullet$\hss}}%
    \graphtemp=.5ex\advance\graphtemp by 0.291in
    \rlap{\kern 3.720in\lower\graphtemp\hbox to 0pt{\hss $\bullet$\hss}}%
    \graphtemp=.5ex\advance\graphtemp by 0.291in
    \rlap{\kern 4.136in\lower\graphtemp\hbox to 0pt{\hss $\bullet$\hss}}%
    \graphtemp=.5ex\advance\graphtemp by 0.707in
    \rlap{\kern 2.965in\lower\graphtemp\hbox to 0pt{\hss $\bullet$\hss}}%
    \graphtemp=.5ex\advance\graphtemp by 0.707in
    \rlap{\kern 3.227in\lower\graphtemp\hbox to 0pt{\hss $\bullet$\hss}}%
    \graphtemp=.5ex\advance\graphtemp by 0.707in
    \rlap{\kern 3.381in\lower\graphtemp\hbox to 0pt{\hss $\bullet$\hss}}%
    \graphtemp=.5ex\advance\graphtemp by 0.707in
    \rlap{\kern 3.643in\lower\graphtemp\hbox to 0pt{\hss $\bullet$\hss}}%
    \special{pa 2965 707}%
    \special{pa 3096 291}%
    \special{pa 3227 707}%
    \special{fp}%
    \special{pa 3381 707}%
    \special{pa 3512 291}%
    \special{pa 3643 707}%
    \special{fp}%
    \special{pa 3096 291}%
    \special{pa 3173 291}%
    \special{fp}%
    \special{pa 3435 291}%
    \special{pa 3797 291}%
    \special{fp}%
    \special{pa 4059 291}%
    \special{pa 4136 291}%
    \special{fp}%
    \graphtemp=.5ex\advance\graphtemp by 0.728in
    \rlap{\kern 3.096in\lower\graphtemp\hbox to 0pt{\hss $\cdots$\hss}}%
    \graphtemp=.5ex\advance\graphtemp by 0.728in
    \rlap{\kern 3.512in\lower\graphtemp\hbox to 0pt{\hss $\cdots$\hss}}%
    \graphtemp=.5ex\advance\graphtemp by 0.832in
    \rlap{\kern 3.096in\lower\graphtemp\hbox to 0pt{\hss $p-1$\hss}}%
    \graphtemp=.5ex\advance\graphtemp by 0.832in
    \rlap{\kern 3.512in\lower\graphtemp\hbox to 0pt{\hss $p-1$\hss}}%
    \graphtemp=.5ex\advance\graphtemp by 0.083in
    \rlap{\kern 3.304in\lower\graphtemp\hbox to 0pt{\hss $t_n$\hss}}%
    \graphtemp=.5ex\advance\graphtemp by 0.104in
    \rlap{\kern 3.928in\lower\graphtemp\hbox to 0pt{\hss $D'$\hss}}%
    \graphtemp=.5ex\advance\graphtemp by 0.312in
    \rlap{\kern 3.304in\lower\graphtemp\hbox to 0pt{\hss $\cdots$\hss}}%
    \graphtemp=.5ex\advance\graphtemp by 0.312in
    \rlap{\kern 3.928in\lower\graphtemp\hbox to 0pt{\hss $\cdots$\hss}}%
    \special{pa 3096 187}%
    \special{pa 3096 83}%
    \special{pa 3137 83}%
    \special{fp}%
    \special{pa 3512 187}%
    \special{pa 3512 83}%
    \special{pa 3449 83}%
    \special{fp}%
    \special{pa 3720 187}%
    \special{pa 3720 83}%
    \special{pa 3761 83}%
    \special{fp}%
    \special{pa 4136 187}%
    \special{pa 4136 83}%
    \special{pa 4074 83}%
    \special{fp}%
    \graphtemp=.5ex\advance\graphtemp by 0.312in
    \rlap{\kern 2.888in\lower\graphtemp\hbox to 0pt{\hss $\cdots$\hss}}%
    \graphtemp=.5ex\advance\graphtemp by 0.499in
    \rlap{\kern 2.472in\lower\graphtemp\hbox to 0pt{\hss $n$ segments\hss}}%
    \graphtemp=.5ex\advance\graphtemp by 0.707in
    \rlap{\kern 2.472in\lower\graphtemp\hbox to 0pt{\hss $C'$\hss}}%
    \graphtemp=.5ex\advance\graphtemp by 1.331in
    \rlap{\kern 4.344in\lower\graphtemp\hbox to 0pt{\hss $\bullet$\hss}}%
    \graphtemp=.5ex\advance\graphtemp by 0.499in
    \rlap{\kern 4.344in\lower\graphtemp\hbox to 0pt{\hss $\bullet$\hss}}%
    \graphtemp=.5ex\advance\graphtemp by 0.499in
    \rlap{\kern 4.344in\lower\graphtemp\hbox to 0pt{\hss $\bullet$\hss}}%
    \special{pa 4136 291}%
    \special{pa 4344 499}%
    \special{pa 4344 603}%
    \special{fp}%
    \special{pa 4136 1540}%
    \special{pa 4344 1331}%
    \special{pa 4344 1227}%
    \special{fp}%
    \special{sh 0.200}%
    \special{ar 4344 915 156 343 0 6.28319}%
    \graphtemp=.5ex\advance\graphtemp by 1.331in
    \rlap{\kern 4.427in\lower\graphtemp\hbox to 0pt{\hss $a$\hss}}%
    \graphtemp=.5ex\advance\graphtemp by 0.499in
    \rlap{\kern 4.427in\lower\graphtemp\hbox to 0pt{\hss $z$\hss}}%
    \graphtemp=.5ex\advance\graphtemp by 0.915in
    \rlap{\kern 4.344in\lower\graphtemp\hbox to 0pt{\hss $R_b$\hss}}%
    \hbox{\vrule depth2.247in width0pt height 0pt}%
    \kern 4.500in
  }%
}%
}
\sp
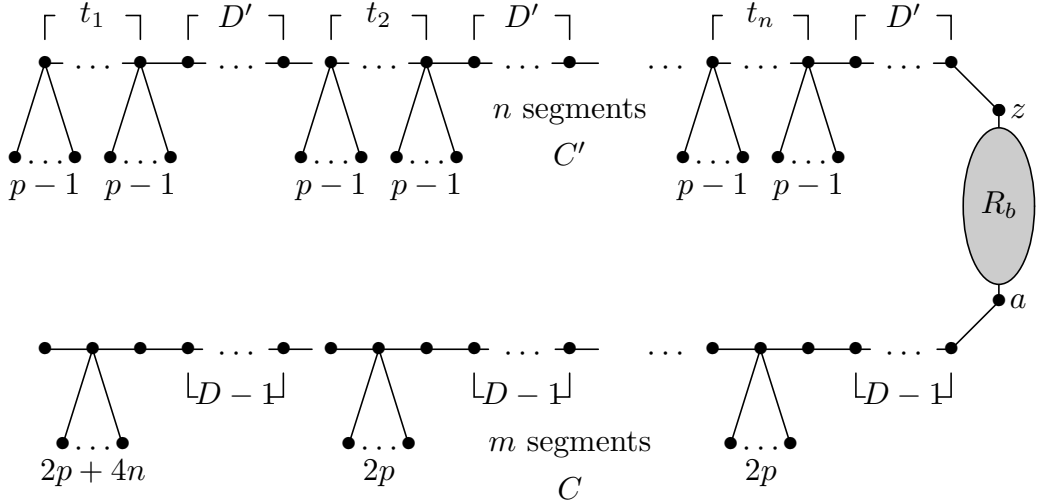
\ce{Fig.~8.  The bug $G$ in the transformation}
\sp

{\bf Schedule yields numbering.}
Suppose that $T$ is solvable; we construct a numbering $f$ of $G$ such that
$B(f)=b$.  As in the proof of Theorem 3, we let $J_i$ be the integer interval
$[ib+1,(i+1)b-1]$.  We first place the vertices from the spine of $C$ in the
positions $\{ib\st 0\le i< \lmb\}$.  For such a vertex in position $ib$, we
place half its leaf neighbors in $J_{i-1}$ and half in $J_i$, using the lowest
$p$ positions in each interval.  For the reflector $R_b$, we use an optimal
numbering in positions $\lmb b,\dots,\lmb b+5b$, with $f(a)=\lmb b$ and
$f(z)=\lmb b+1$ (see Fig.~7).

Since $T$ is solvable, there exist index sets $\{I_j\st 1\le j\le m\}$ assigning
jobs to processors such that $\UE j1m I_j = \{1,\dots,n\}$ and that
$\SM i{I_j} t_i\le D$ for each $j$.  We may increase the task execution times to
obtain $\SM i{I_j} t_i= D$ for each $j$, because this does not change $b$ and
the new corresponding bug contains the original bug as an induced subgraph.
The segments of $C'$ indexed by $I_j$ thus have exactly $(\SM i{I_j} t_i)p = Dp$
vertices consisting of leaves and their neighbors; let $S_j$ denote this set
of vertices in $C'$.

We have already placed the vertices of the $j$th segment of $C$ in some of the
positions from $(j-1)(D+2) b$ to $(j(D+2)-1)b$, using all the multiples of $b$
and some positions in $J_{(j-1)(D+2)}$ and $J_{(j-1)(D+2)+1}$.  We now place the
vertices of $S_j$ in $J_{(j-1)(D+2)+2},\dots,J_{j(D+2)-1}$, with $J_i$ receiving
one spine vertex at position $ib+1$, followed immediately by its $p-1$ leaf
neighbors in positions $ib+2,\dots,ib+p$.

In each such $J_i$, we have assigned $p$ positions and still have $b-1-p=2n$
unassigned.  This also holds for $J_{j(D+2)}$ and $J_{j(D+2)+1}$ such that
$1\le j<m$, where we have placed leaves from $C$ in the lowest positions (we
have filled $J_0$ and $J_1$ completely).  Thus there remain $2n(\lmb-2)$
positions unassigned.

To complete the numbering, we must assign positions to the remaining vertices
from the spine of $C'$.  The remaining vertices consist of $n$ paths, each of
order $D'=2(\lmb-2)$.  We will place two vertices from each path into each $J_i$
for $2\le i<\lmb$.  For fixed $k\in\{1,\ldots,n\}$, we place the $k$th path
into $L\cup U$, where $L = \{ib+p+2k-1\st 2\le i<\lmb\}$ and 
$U = \{ib+p+2k\st 2\le i<\lmb\}$.

We place the two endpoints of the path in the intervals $J,J'$ containing their
neighbors on the spine of $C'$.  From the higher desired interval $J'$, the path
moves up to $J_{\lmb-1}$ via $L$.  It then switches to $U$ and moves down to
$J'$ via $U$.  From there down to its entrance to $J$, it uses the positions of
$U$ and $L$ in each interval.  It moves from $J$ down to $J_2$ via $U$,
switches to $L$ in $J_2$, and moves back up to $J$ via $L$, where
it ends (see Fig.~9).  Because successive positions within $L$ or within $U$
differ by exactly $b$, we have completed a numbering showing that $B(G)\le b$.

\sp
\gpic{
\expandafter\ifx\csname graph\endcsname\relax \csname newbox\endcsname\graph\fi
\expandafter\ifx\csname graphtemp\endcsname\relax \csname newdimen\endcsname\graphtemp\fi
\setbox\graph=\vtop{\vskip 0pt\hbox{%
    \graphtemp=.5ex\advance\graphtemp by 0.143in
    \rlap{\kern 0.323in\lower\graphtemp\hbox to 0pt{\hss $\bullet$\hss}}%
    \graphtemp=.5ex\advance\graphtemp by 0.143in
    \rlap{\kern 0.682in\lower\graphtemp\hbox to 0pt{\hss $\bullet$\hss}}%
    \graphtemp=.5ex\advance\graphtemp by 0.143in
    \rlap{\kern 1.040in\lower\graphtemp\hbox to 0pt{\hss $\bullet$\hss}}%
    \graphtemp=.5ex\advance\graphtemp by 0.143in
    \rlap{\kern 1.399in\lower\graphtemp\hbox to 0pt{\hss $\bullet$\hss}}%
    \graphtemp=.5ex\advance\graphtemp by 0.143in
    \rlap{\kern 1.758in\lower\graphtemp\hbox to 0pt{\hss $\bullet$\hss}}%
    \graphtemp=.5ex\advance\graphtemp by 0.143in
    \rlap{\kern 2.117in\lower\graphtemp\hbox to 0pt{\hss $\bullet$\hss}}%
    \graphtemp=.5ex\advance\graphtemp by 0.143in
    \rlap{\kern 2.475in\lower\graphtemp\hbox to 0pt{\hss $\bullet$\hss}}%
    \graphtemp=.5ex\advance\graphtemp by 0.143in
    \rlap{\kern 2.834in\lower\graphtemp\hbox to 0pt{\hss $\bullet$\hss}}%
    \graphtemp=.5ex\advance\graphtemp by 0.143in
    \rlap{\kern 3.193in\lower\graphtemp\hbox to 0pt{\hss $\bullet$\hss}}%
    \graphtemp=.5ex\advance\graphtemp by 0.143in
    \rlap{\kern 3.552in\lower\graphtemp\hbox to 0pt{\hss $\bullet$\hss}}%
    \graphtemp=.5ex\advance\graphtemp by 0.143in
    \rlap{\kern 3.910in\lower\graphtemp\hbox to 0pt{\hss $\bullet$\hss}}%
    \graphtemp=.5ex\advance\graphtemp by 0.000in
    \rlap{\kern 0.502in\lower\graphtemp\hbox to 0pt{\hss $J_2$\hss}}%
    \graphtemp=.5ex\advance\graphtemp by 0.000in
    \rlap{\kern 1.578in\lower\graphtemp\hbox to 0pt{\hss $J$\hss}}%
    \graphtemp=.5ex\advance\graphtemp by 0.000in
    \rlap{\kern 2.655in\lower\graphtemp\hbox to 0pt{\hss $J'$\hss}}%
    \graphtemp=.5ex\advance\graphtemp by 0.000in
    \rlap{\kern 3.731in\lower\graphtemp\hbox to 0pt{\hss $J_{\lmb -1}$\hss}}%
    \graphtemp=.5ex\advance\graphtemp by 0.359in
    \rlap{\kern 0.000in\lower\graphtemp\hbox to 0pt{\hss $U$\hss}}%
    \graphtemp=.5ex\advance\graphtemp by 0.502in
    \rlap{\kern 0.000in\lower\graphtemp\hbox to 0pt{\hss $L$\hss}}%
    \special{pn 8}%
    \special{pa 466 143}%
    \special{pa 538 143}%
    \special{fp}%
    \special{pa 1901 143}%
    \special{pa 1973 143}%
    \special{fp}%
    \special{pa 2260 143}%
    \special{pa 2332 143}%
    \special{fp}%
    \special{pa 3695 143}%
    \special{pa 3767 143}%
    \special{fp}%
    \special{ar 646 -80 287 287 0.895665 2.245928}%
    \special{ar 1004 -80 287 287 0.895665 2.245928}%
    \special{ar 1363 -80 287 287 0.895665 2.245928}%
    \special{ar 717 -80 287 287 0.895665 2.245928}%
    \special{ar 1076 -80 287 287 0.895665 2.245928}%
    \special{ar 1435 -80 287 287 0.895665 2.245928}%
    \special{ar 1758 -104 287 287 1.047198 2.094395}%
    \special{ar 2117 -104 287 287 1.047198 2.094395}%
    \special{ar 2511 -80 287 287 0.895665 2.245928}%
    \special{ar 2870 -80 287 287 0.895665 2.245928}%
    \special{ar 3229 -80 287 287 0.895665 2.245928}%
    \special{ar 3587 -80 287 287 0.895665 2.245928}%
    \special{ar 2798 -80 287 287 0.895665 2.245928}%
    \special{ar 3157 -80 287 287 0.895665 2.245928}%
    \special{ar 3516 -80 287 287 0.895665 2.245928}%
    \special{pa 2619 502}%
    \special{pa 2978 502}%
    \special{pa 3336 502}%
    \special{pa 3695 502}%
    \special{pa 3767 359}%
    \special{pa 3408 359}%
    \special{pa 3049 359}%
    \special{fp}%
    \special{pa 3049 359}%
    \special{pa 2691 359}%
    \special{pa 2332 359}%
    \special{pa 2260 502}%
    \special{pa 1973 359}%
    \special{pa 1901 502}%
    \special{pa 1614 359}%
    \special{fp}%
    \special{pa 1614 359}%
    \special{pa 1256 359}%
    \special{pa 897 359}%
    \special{pa 538 359}%
    \special{pa 466 502}%
    \special{pa 825 502}%
    \special{pa 1184 502}%
    \special{fp}%
    \special{pa 1184 502}%
    \special{pa 1543 502}%
    \special{fp}%
    \hbox{\vrule depth0.646in width0pt height 0pt}%
    \kern 4.000in
  }%
}%
}
\sp
\ce{Fig.~9.  Numbering a path}
\sp

{\bf Numbering yields schedule.}
Conversely, suppose that $G$ has a numbering $f$ with $B(f)=b$ (there is
no better numbering, since $B(R_b)=b$.  We prove that an optimal numbering
of $G$ must have essentially the form described above, from which we obtain
a positive solution for $T$.

The positions within $b$ of the vertex $x$ of degree $2+2p+4n=2b$ on the spine
of $C$ must be occupied by the neighbors of $x$, and no edge can stretch across
this interval.  Thus $x$ and these vertices occupy positions at one end of the
numbering; by symmetry, we may assume that these positions are $0,\dots,2b$.
Similarly, the reflector $R_b$ must occupy $5b+1$ consecutive positions, with
its peripheral vertices within $b$ of the end, and no edge can stretch across
these positions.  Thus we may assume that $R-b$ occupies positions
$\lmb b,\dots,(\lmb+5)b$ and that $\{f(a),f(z)\}=\{\lmb b,\lmb b+1\}$.

Among the remaining positions, which must be filled since $G$ has
$(\lmb+5)b+1$ vertices, there is a path of length $\lmb$ from a leaf
neighbor of $x$ to $a$.  Since there is a leaf neighbor of $x$ in position 0
and $a$ has position at least $\lmb b$, the vertices of this path must occupy
the positions $\{ib\st0\le i\le \lmb\}$.

It remains to assign tasks to sets $I_j$ such that $\SM i{I_j} t_i \le D$.
Let $z_j=[1+j(D+2)]b$ for $0\le j\le m$.  We say that task $i$ belongs to
set $I_j$ if and only if $z_{j-1}<f(v)<z_j$ for some non-leaf task vertex $v$ in
the $i$th segment of $C'$.  We show first that task $i$ belongs to only one set;
suppose not.  Since the non-leaf task vertices in the $i$th segment of $C'$
induce a connected subgraph (a path), there must exist adjacent non-leaf task
vertices $u,v$ such that $f(u)<z_j<f(v)$ for some $j$.  Since $u$ and $v$ are
adjacent, we have $f(v)-f(u)\le b$.  Now $u,v,f^{-1}(z_j)$ and their neighbors
all have positions in the interval $[f(u)-b,f(v)+b]$.  There are $4p+3$ of
these vertices, but at most $3b+1$ positions in the interval.  From $b=p+1+2n$
and $p>8n$ we obtain $4p+3>3b+1$, and we cannot place $4p+3$ vertices into
$3b+1$ positions.

The positions outside $[z_0,z_m]$ are already filled, so we have assigned each
task to exactly one set $I_j$.  We must show that $\SM j{I_j} t_i\le D$ for
each $j$.  Into the interval $[z_{j-1}-b,z_j+b]$, we have now placed all the
task vertices for tasks in $I_j$, $D+5$ vertices from the spine of $C$,
and $4p$ leaves of $C$ (when $j=m$, some of this count is replaced by vertices
of $R_b$).  The number of vertices is $(\SM j{I_j} t_i)p+D+5+4p$, and the
number of positions in the interval is $(D+2)b+1+2b$.  Thus
$\SM j{I_j} t_i \le [(D+4)(p+1+2n)-(D+5+4p)]/p = D+[2n(D+4)-1]/p$.
Since $p>2n(D+4)$, we conclude that $\SM j{I_j} t_i \le D$.  \qed

The paradigm in the proof of Theorem 4 applies more generally.
All we need is a graph to play the role of the reflector.  This graph $R_b'$
will have bandwidth $b$ and two special vertices such that every optimal
numbering puts those two vertices in positions near one end.
We can then use $R_b'$ in place of $R_b$ to form a bug-like graph and
follow the proof of Theorem 4.

For example, we use this approach to prove that bandwidth is NP-complete on
a class of trees that are tolerance graphs.  A {\it near-caterpillar}
is a tree having a single path that includes all but one of the non-leaf
vertices.  For our reflector with bandwidth $b$ we use a near-caterpillar with
$4b+1$ vertices very similar to the near-caterpillar $T_b$ of Fig.~3.
Define $R_b'$ for $b$ even to be the same as $T_b$ in Fig.~3 except that
the $4b-3$ leaves are redistributed among the sets $X,Y,Z,W$ so that the sizes
of $X,Y,Z,W$ are $b/2,b,b/2,2b-3$, respectively.  The bandwidth of $T_b$ exceeds
$b$; the bandwidth of $R_b'$ is its local density $b$, but this is achievable
only by putting specified vertices near one end.

\LM 5.
In every optimal numbering $f$ of the near-caterpillar $R_b'$ using positions
$0,\dots,4b$, the vertices of $X\cup Z$ are all below $b$ or all above $3b$.
\PF
The local density of $R_b'$ is $b$, and we have a numbering achieving bandwidth
$b$ in which $x,z,w,y$, are numbered $b$, $b+1$, $2b$, $3b$, respectively.  Now
consider an arbitrary optimal numbering $f$.  By symmetry, we may assume that
$f(w)<f(y)$.  Since $w$ has degree $2b$, its neighbors fill the positions within
$b$ to each side of $f(w)$.  Hence all of $Y$ is outside this interval, and no
edges stretch across, so $f(y)= f(w)+b$ and $Y$ occupies $[f(w)+b+1,f(w)+2b]$.
Since no edge involving $x$ or $z$ can stretch across this interval, all
of $X\cup Z$ must be on the other end, below $f(w)-b$.  We thus have $2b$
vertices to each side of $f(w)$, and all of $X\cup Z$ is below $b$.  \qed

By following the argument of Theorem 4, bandwidth is NP-complete for
near-caterpillars.  It is immediate that every near-caterpillar
is a tolerance graph, since a tree is a tolerance graph if and only if it does
not contain the tree obtained from the claw $K_{1,3}$ by subdividing each
edge twice [7].  This tree is forbidden from near-caterpillars, since every
path in it misses at least two non-leaf vertices.  Thus bandwidth is NP-complete
for a subclass of tolerance graphs.

\SH
{References}
\frenchspacing
\BP [1]
S.F. Assmann, G.W. Peck, M.M. Sys\l o, and J. Zak, The bandwidth of
caterpillars with hairs of length $1$ and $2$.  \SIAD\ 2(1981), 387--393. 
\BP [2]
P.Z.~Chinn, J.~Chv\'atalov\'a, A.K.~Dewdney, and N.E.~Gibbs,
The bandwidth problem for graphs and matrices - A survey.
\JGT\ 6(1982), 223--254.
\BP [3]
F.R.K.~Chung, Labelings of graphs,
in {\it Selected Topics in Graph Theory, III} (L.~Beineke and R.~Wilson, eds.).
(Academic Press 1988), 151--168.
\BP [4]
M.R.~Garey, R.L.~Graham, D.S.~Johnson, and D.E.~Knuth, Complexity results for
bandwidth minimization.  \SIAP\ 34(1978), 477--495.
\BP [5]
M.R.~Garey and D.S.~Johnson, {\it Computers and Intractability: A Guide to
the Theory of NP-Completeness}, (W.H.~Freeman, 1979).
\BP [6]
M.C.~Golumbic and C.L.~Monma, A generalization of interval graphs with
tolerances, {\it Proc. 13th SE Conf. Comb. G.T. Comp.}, \CNum\ 35(1982),
37--46.
\BP [7]
M.C.~Golumbic, C.L.~Monma, and W.T.~Trotter, Tolerance graphs,
\DAM\ 9(1984), 157--170.
\BP [8]
F.~Harary, A characterization of block-graphs,
{\it Canad.\ Math.\ Bull.} 6(1963), 1--6.
\BP [9]
D.J.~Kleitman and R.V.~Vohra, Computing the bandwidth of interval graphs,
\SIDM\ 3(1990), 373--375.
\BP [10]
F.~Makedon and I.H.~Sudborough, On minimizing width in linear layouts,
\DAM\ 23(1989), 243--265.
\BP [11]
Z.~Miller, The bandwidth of caterpillar graphs, {\it Proc.\ Southeastern
Conf.{}}, {\it Congressus Numerantium} 33(1981), 235--252.
\BP [12]
B.~Monien, The bandwidth minimization problem for caterpillars with hair
length 3 is NP-complete, \SIAD\ 7(1986), 505--512.
\BP [13]
D.O.~Muradyan, A polynomial algorithm for finding the bandwidth of interval
graphs, {\it Dokl. AN Arm. SSR} 82(1986), 64--66.
\BP [14]
C.H.~Papadimitriou, the NP-completeness of the bandwidth minimization problem.
{\it Computing} 16(1976), 263--270.
\BP [15]
L. Smithline, Bandwidth of the complete $k$-ary tree,
\DM\ 142(1995), 203--212.
\BP [16]
A.P.~Sprague, An $O(n\log n)$ algorithm for bandwidth of interval graphs,
\SIDM\  {7}(1994), 213--220. 
\BP [17]
M.M.~Sys\l o and J.~Zak, The bandwidth problem: critical subgraphs and the
solution for caterpillars, \ADM\ 16(1982), 281--286.  See also Comp.\ Sci.\
Dept.\ Report CS-80-065, Washington State Univ.\ (1980).
\bye